\documentclass[11pt]{article}

\usepackage[a4paper,margin=1in]{geometry}
\usepackage{amsmath,amssymb,amsthm,mathtools}
\usepackage{bm}
\usepackage{microtype}
\usepackage{enumitem}
\usepackage{needspace}
\usepackage{pgfplots}
\pgfplotsset{compat=1.18}
\usepackage[colorlinks=true,linkcolor=blue,citecolor=blue,urlcolor=blue]{hyperref}

\newcounter{proofserial}
\newcounter{proofstep}[proofserial]

\newcommand{\stepheading}[1]{%
  \refstepcounter{proofstep}%
  \emph{Step \theproofstep: #1}\space}
\newcommand{\proofstep}[1]{\par\indent\stepheading{#1}}

\newlength{\proofheadingindent}
\makeatletter
\renewenvironment{proof}[1][\proofname]{\par
  \stepcounter{proofserial}%
  \pushQED{\qed}%
  \normalfont \topsep6\p@\@plus6\p@\relax
  \trivlist
  \item[\hskip\labelsep\hspace{\proofheadingindent}%
    \itshape #1\@addpunct{.}]\ignorespaces
}{\popQED\endtrivlist\@endpefalse}
\makeatother

\numberwithin{equation}{section}

\newtheorem{theorem}{Theorem}[section]
\newtheorem{proposition}[theorem]{Proposition}
\newtheorem{lemma}[theorem]{Lemma}
\newtheorem{corollary}[theorem]{Corollary}
\theoremstyle{definition}
\newtheorem{assumption}{Assumption}[section]
\theoremstyle{definition}
\newtheorem{remark}[theorem]{Remark}
\newtheorem{example}[theorem]{Example}

\newcommand{\R}{\mathbb{R}}

\newcommand{\Prob}{\mathbb{P}}          
\newcommand{\E}{\mathbb{E}}
\newcommand{\F}{\mathcal{F}}
\newcommand{\Eset}{\mathcal{E}}
\newcommand{\Gset}{\mathcal{G}}
\newcommand{\Hset}{\mathcal{H}}
\newcommand{\Bset}{\mathcal{B}}
\newcommand{\one}{\mathbf{1}}
\newcommand{\dd}{\,\mathrm{d}}

\newcommand{\cN}{\mathcal{N}}
\newcommand{\cL}{\mathcal{L}}
\newcommand{\cR}{\mathcal{R}}
\newcommand{\dist}{\operatorname{dist}}
\newcommand{\Var}{\operatorname{Var}}

\newcommand{\argmin}{\operatorname*{arg\,min}}
\newcommand{\geom}{\mathrm{geom}}

\newcommand{\Op}{\mathrm{op}}

\title{Terminal Localization of Empirical Reverse Diffusions\\ under Singular Score Errors%
\thanks{This work was supported by the National Natural Science
  Foundation of China under Grant No.~12401586. The work was completed during the
  author's visit to Friedrich-Alexander-Universit\"at Erlangen-N\"urnberg, and the
  author thanks Professor Enrique Zuazua for constructive comments on the paper.}}
\author{Zhonghua Liao\thanks{School of Mathematics, Sichuan University, Chengdu, China.
  \texttt{zhonghualiao@yeah.net}\\
  \emph{AI-usage disclosure.} This paper was written with extensive use of
  ChatGPT, including but not limited to the proofs of the original paper and their
  improvement, the language editing, the numerical plots, and the search for
  references. The author's contribution consists mainly of the human processes that
  AI cannot currently replace: communicating with ChatGPT to propose and refine the
  research directions, guiding it through the proofs, reading, checking, and
  improving the drafts it produced, and verifying that the cited references are
  appropriate.}}
\date{}

\begin{document}
\maketitle
\begin{abstract}
For a finite empirical target, the exact reverse diffusion is driven toward
training points by an attraction whose strength diverges as the forward
noising time decreases to zero. We study whether this terminal localization
persists under progressively measurable score errors that may themselves be
singular. For reverse Ornstein--Uhlenbeck diffusions, from an arbitrary
initial law and under a deterministic error envelope \(E\) imposed along the
perturbed trajectory satisfying
\(t\int_t^T E(\sigma)/\sigma\,d\sigma\to0\), we prove almost-sure
convergence to a training point, together with explicit uniform bounds on a
terminal window with probability arbitrarily close to one. Every training
point has positive terminal probability, although the terminal selection
probabilities need not match the empirical weights. The results cover
power-law errors \(E(t)=O(t^{-\gamma})\), \(0<\gamma<1\), including infinite
squared-error energy. For suitable errors in this range, the exact and
perturbed full path laws can be mutually singular while their terminal laws
share the same finite support. For suitable switched perturbations at the
critical scale \(t^{-1}\), arbitrarily small positive amplitudes can produce
a nonzero terminal displacement with arbitrarily high probability. These
results quantify stochastic robustness of empirical terminal attraction
without requiring absolute continuity of the full path laws.
\end{abstract}

\section{Introduction}\label{sec:intro}

Diffusion models have become a central framework for generative modeling,
with successful applications to image synthesis \cite{rombach2022},
audio generation \cite{kong2021}, and three-dimensional molecule generation
\cite{hoogeboom2022}.
They generate data by reversing a noise-adding diffusion, with a drift involving the score, the logarithmic gradient of the forward reference density
\cite{anderson1982,cattiaux2023,haussmann1986,ho2020,song2019,song2021}.
For a finite empirical target, the exact continuous-time reverse diffusion converges to one of the training points as the forward noising time \(t\) decreases to zero during generation.
This attraction to training samples is often called memorization in the diffusion-model literature
\cite{baptista2025,biroli2024,li2024goodscore}; in this paper the property we prove is called
terminal localization. It concerns the endpoint of a single trajectory, and we do not use
memorization as a synonym for faithful sampling.
In practical applications, the exact score is generally unavailable, so the
reverse dynamics use an approximate score learned from data
\cite{song2019,song2021}.
This approximation changes the reverse drift, and the resulting error may
become singular near the generation endpoint \(t=0\).
We ask whether terminal localization persists under such perturbations and how the behavior depends on the rate at which the score error diverges.

A natural approach to this question is Girsanov's theorem
\cite{girsanov1960}, which compares the exact and perturbed diffusions by
changing measure on path space. This approach is used in early quantitative
convergence analyses of diffusion models
\cite{chen2023sampling,debortoli2021}. Pidstrigach
\cite[Assumption~2 and Theorem~2]{pidstrigach2022} also uses Girsanov weights
to establish terminal limits and equivalence of terminal laws. However, the
stochastic exponential defining the change of measure must be a true
martingale, and extending this comparison to the terminal time requires
uniform integrability of the density process. Infinite-energy perturbations can instead produce mutually
singular full path laws, even when Girsanov's theorem applies on every interval
strictly before the endpoint. This leaves open whether terminal localization
persists when the exact and perturbed full path laws are mutually singular.

To address this question, we examine the singular attraction generated by
the empirical score. Its increasing strength near the terminal time suggests
measuring score errors through their contraction-weighted accumulated
response. Since attraction toward an individual training point is effective
only near that point, a global argument is also needed to bring the
trajectory into a neighborhood where this contraction can maintain
localization.

For persistent power-law envelopes \(E(t)\asymp t^{-\gamma}\), where \(t\)
is the forward noising time and generation runs toward \(t=0\), the results distinguish the following regimes:
\begin{center}
\begin{tabular}{@{}p{0.30\linewidth}@{\hspace{0.035\linewidth}}p{0.30\linewidth}@{\hspace{0.035\linewidth}}p{0.30\linewidth}@{}}
\multicolumn{1}{c}{\(\gamma<1/2\)} &
\multicolumn{1}{c}{\(1/2\le\gamma<1\)} &
\multicolumn{1}{c}{\(\gamma=1\)}\\[3pt]
{\raggedright Finite energy; path-law comparison by Girsanov's theorem\par} &
{\raggedright Localization holds; path laws can be mutually singular\par} &
{\raggedright Localization can fail for arbitrarily small positive amplitudes\par}
\end{tabular}
\end{center}
The middle regime admits deterministic infinite-energy errors with mutually
singular path laws (Lemma~\ref{lem:nonui}); the critical regime concerns the
switched perturbations of Proposition~\ref{prop:critical}, whose construction
also allows arbitrary initialization laws.
Figure~\ref{fig:response-energy-illustration} and Table~\ref{tab:tstop}
illustrate these three regimes by direct simulation of \eqref{eq:main-sde}.

We next introduce the main model studied in this paper and state the main theorem.

Let
$S\triangleq \{y_1,\ldots,y_N\}\subset\R^d$ be a finite set of distinct points, with
$N\ge2$ and $d\ge1$, and let
$u_0\triangleq N^{-1}\sum_{i=1}^N\delta_{y_i}$ be its empirical data measure,
so $S=\operatorname{supp}(u_0)$. Its Ornstein--Uhlenbeck smoothing has density
\begin{equation}\label{eq:model}
\begin{aligned}
 w(x,t)&\triangleq\frac1N\sum_{i=1}^N\varphi_{h_t}(x-a_ty_i),\qquad t>0,\\
 a_t&\triangleq e^{-t},\qquad h_t\triangleq1-e^{-2t}.
\end{aligned}
\end{equation}
Here $\varphi_h$ denotes the centered Gaussian density in $\R^d$ with covariance
$hI_d$. We use $w(x,t)$ only for the density and $w(t)\triangleq w(x,t)\dd x$ only
for the corresponding probability measure; the symbol $w$ always denotes this OU
reference.
The initial condition $w(\cdot,0)=u_0$ is understood in the sense of probability
measures. Equivalently, $w(t)$ is the law of the forward noising process
$\dd X_t^{\mathrm f}=-X_t^{\mathrm f}\dd t+\sqrt2\,\dd B_t^{\mathrm f}$,
$X_0^{\mathrm f}\sim u_0$.
Its exact empirical score is
\begin{equation}\label{eq:score}
\begin{aligned}
 s(x,t)&\triangleq\nabla_x\log w(x,t)=\frac{\mu(x,t)-x}{h_t},\\
 \mu(x,t)&\triangleq\E[a_tX_0^{\mathrm f}\mid X_t^{\mathrm f}=x]
 =\sum_{i=1}^N\omega_i(x,t)\,a_ty_i,\\
 \omega_i(x,t)&\triangleq
 \frac{\exp\bigl(-|x-a_ty_i|^2/(2h_t)\bigr)}
 {\sum_{j=1}^N\exp\bigl(-|x-a_ty_j|^2/(2h_t)\bigr)}.
\end{aligned}
\end{equation}
The weights $\omega_i$ are the conditional probabilities of the mixture
components, and $\mu(x,t)$ is the conditional mean of the moving centers
$a_ty_i$ given the noisy state. Throughout, $s$ denotes this exact score, and a
learned score is denoted $s_\theta$. We write
\[
\Lambda\triangleq\min_{i\ne j}|y_i-y_j|>0,
\qquad R\triangleq\max_{1\le i\le N}|y_i|.
\]
Since $h_t\sim2t$, the coefficient $h_t^{-1}$ in \eqref{eq:score} has order
$t^{-1}$ near the generation endpoint.

Fix $T>0$. We use $t\in(0,T]$ for forward noising time, with generation
running from $t=T$ toward $t=0$. For stochastic analysis, set $r=T-t$,
the increasing generation time. Let $(\Omega,\F,(\F_r)_{0\le r<T},\Prob)$
be a filtered probability space with a complete and right-continuous
filtration, and let $B=(B_r)_{0\le r<T}$ be a $d$-dimensional
$(\F_r)$-Brownian motion. Set
$\F_{T-}\triangleq\sigma(\bigcup_{0\le r<T}\F_r)$. We define generation
through the increasing-time It\^o SDE
\begin{equation}\label{eq:main-sde}
\begin{aligned}
 \dd\widehat X_r
 &=\bigl[\widehat X_r+2s(\widehat X_r,T-r)+2e_r\bigr]\dd r
   +\sqrt2\,\dd B_r,\\
 &\qquad 0\le r<T,\qquad \widehat X_0\sim v_T,
\end{aligned}
\end{equation}
where $\widehat X_0$ is $\F_0$-measurable and $v_T$ is an arbitrary
initialization probability measure on $\R^d$. Define
\[
 X_t\triangleq\widehat X_{T-t},\qquad
 v(t)\triangleq\cL(X_t),\qquad 0<t\le T.
\]
Thus $X_T=\widehat X_0\sim v_T$. Throughout, $v(t)=\cL(X_t)$ is always a
probability measure, namely the generated law; it need not have a density, and we
never write $v(t)$ for a density. The drift uses the prescribed score
$s=\nabla\log w$ of the forward OU density. When the error induced by a learned
score is Markovian, it has the form
$e_r=s_\theta(\widehat X_r,T-r)-s(\widehat X_r,T-r)$; the analysis uses only the
progressively measurable error $e$ of Assumption~\ref{ass:localization}, which also
covers path-dependent, non-Markov perturbations.

\begin{center}
\fbox{\parbox{0.92\linewidth}{\small
\textbf{Convention on time.} The forward noising time $t\in(0,T]$ \emph{decreases}
to $0$ during generation, while the generation time $r=T-t$ \emph{increases} to $T$.
The two processes are the same path, $X_t=\widehat X_{T-t}$, so that $t\downarrow0$
and $r\uparrow T$ refer to the same instant, the generation endpoint. A terminal
window is $0<t\le b$, equivalently $T-b\le r<T$, with $0<b<T$. The word
\emph{terminal} always refers to $t\downarrow0$, while the initialization law
$v_T$ is the law at the \emph{start} of generation, $t=T$. Terminal localization
means that the trajectory stays near the same training sample throughout this final
generation stage.}}
\end{center}

For terminal localization, we impose the following assumption:
\begin{assumption}[Solution and score-error conditions]\label{ass:localization}
\begin{enumerate}[label=\textup{(\roman*)},ref=\roman*,leftmargin=2em]
\item \label{ass:local}
The error \(e\triangleq(e_r)_{0\le r<T}\) is progressively measurable
with respect to \((\F_r)\). There exists a deterministic, nonnegative,
measurable function \(E:(0,T]\to[0,\infty)\) such that
\begin{equation}\label{eq:envelope}
E\in L^2([\delta,T])\quad\text{for every }0<\delta<T,
\qquad |e_r|\le E(T-r)\quad\dd r\otimes\dd\Prob\text{-a.e.}
\end{equation}
For every \(0<\delta<T\), equation~\eqref{eq:main-sde} admits a
continuous adapted solution on \([0,T-\delta]\), and these solutions
agree on overlapping intervals.  
\item \label{ass:terminal}
The envelope \(E\) in part~\textup{(\ref{ass:local})} satisfies
\begin{equation}\label{eq:weighted-vanish}
t\int_t^T\frac{E(\sigma)}{\sigma}\dd \sigma \longrightarrow0
\qquad\text{as }t\downarrow0.
\end{equation}
\end{enumerate}
\end{assumption}
\begin{remark}[Where the bound is imposed]\label{rem:where}
	The envelope in \eqref{eq:envelope} bounds the error along the perturbed trajectories,
	not its $L^2(w(t))$ norm under the reference law. An $L^2(w(t))$ estimate does not imply
	\eqref{eq:envelope}; see Proposition~\ref{prop:global-localization} and
	Corollary~\ref{cor:mean-integrated} for the weaker, realized conditions that suffice
	for qualitative localization.
\end{remark}
We call $E$ the \emph{score-error envelope}; the squared-error energy
$\int E^2$ is a different quantity and is never denoted by $E$.

\begin{theorem}[Terminal localization, simplified form]
\label{thm:intro-localization}
Suppose that Assumption~\ref{ass:localization} holds.
For every \(\alpha\in(0,1)\), there exist a constant
\(b\in(0,T)\), an \(\F_{T-b}\)-measurable index \(I_b\in\{1,\ldots,N\}\), and a deterministic
continuous function \(\cR_{b,\alpha}:(0,b]\to[0,\infty)\) tending to zero as
\(t\downarrow0\), such that, with probability at least \(1-\alpha\),
\[
 \bigl|X_t-a_ty_{I_b}\bigr|\le \cR_{b,\alpha}(t)
 \qquad\text{for every }0<t\le b.
\]
On this event, \(a_ty_{I_b}\) is the unique nearest moving center throughout
the interval, and \(X_t\to y_{I_b}\) as \(t\downarrow0\).
Moreover, \(X_t\) converges almost surely to a training point as
\(t\downarrow0\).
\end{theorem}

\begin{remark}[Role of the weighted error condition]
\label{rem:integrated-response}
The weighted condition~\eqref{eq:weighted-vanish} permits the score error to grow near the terminal
time, even with infinite squared-error energy. For example, it holds for
\(E(t)=O(t^{-\gamma})\) whenever \(0<\gamma<1\). The increasing attraction
toward training points can suppress the displacement caused by an unbounded
score error, and the weighted integral measures its accumulated effect.

Once the trajectory enters a neighborhood of a moving training point, the
proof must ensure that it stays near the same point until the terminal time.
The local propagation estimate bounds the accumulated error by a constant
multiple of the same weighted integral restricted to the terminal interval.
The weighted condition~\eqref{eq:weighted-vanish} makes this contribution uniformly small on a
sufficiently short terminal interval. Together with the
Brownian and geometric bounds, this excludes a first exit from the neighborhood
and maintains localization.
\end{remark}

The explicit terminal estimate and the full-support conclusion are given in
Theorem~\ref{thm:main} and Corollary~\ref{cor:positive-terminal-support},
respectively.

Previous studies have addressed different aspects of terminal behavior in
diffusion models. Pidstrigach
\cite[Assumption~2 and Theorem~2]{pidstrigach2022} proves existence of
terminal limits for perturbed reverse diffusions and equivalence with the
exact terminal law under uniform integrability of the Girsanov weights and
the other regularity hypotheses. De Bortoli
\cite[Theorem~1]{debortoli2022} derives Wasserstein estimates for discretized
reverse diffusions with singular score errors, stopping before the terminal
time. For empirical targets, Baptista et al.
\cite[Theorem~4.3]{baptista2025} study attraction for the empirical-score ODE,
including possible Voronoi-boundary limits. Wan et al.
\cite[Theorem~5.3 and Propositions~5.9--5.10]{wang2025} establish convergence
and absorption for flow-matching ODEs, including results for smooth
approximate denoisers with uniformly vanishing spatial error.

\par Liu and Zuazua \cite{liu2025} analyze generative dynamics from the heat-flow PDE
viewpoint. They study the Fokker--Planck evolution of the generated density, obtain
distributional concentration and exact-score results, and quantify score
approximation through a relative-entropy estimate. Their conclusions therefore
concern the generated \emph{law} and are measured in $W_2$ and in entropy.
The present paper proves a different kind of statement for a closely related
process. After the classical self-similar change of variables
$w(x,t)=e^{dt}u(e^{t}x,(e^{2t}-1)/2)$, where $u$ solves $\partial_s u=\Delta u$
with $u(\cdot,0)=u_0$, the OU reference of \eqref{eq:model} corresponds to the heat
flow of \cite{liu2025}; we keep the OU parametrization because the stochastic
perturbation is defined on paths, not on densities. The object we control is a single realization: we prove
almost-sure convergence of $X_t$ to a training point together with a deterministic
terminal bound that holds uniformly on $(0,b]$ (Theorem~\ref{thm:main}). This
conclusion does not require the perturbed and exact laws to be comparable, and it
applies in the mutually singular regime of Lemma~\ref{lem:nonui}. The two sets of results are thus complementary rather than comparable: \cite{liu2025}
quantifies how the generated law concentrates, while Theorem~\ref{thm:main}
quantifies where an individual trajectory terminates and how fast.

Other related results concern the geometry of singular targets and the
stability of diffusion dynamics. Li et al. \cite[Theorem~4.1]{li2025} study
support concentration under marginal-score discrepancies at the
inverse-variance scale. Chandramoorthy
and de Clercq \cite{chandramoorthy2025} analyze linearized perturbation
responses, and Cao, Wu and Chandramoorthy \cite{cao2026guidance} establish
support convergence for guided SDEs with exact scores under geometric
assumptions. Khelifa, Turner and Venkataramanan \cite{khelifa2026geometry}
study score-error components that preserve marginal dynamics. Further work
concerns data geometry \cite{chen2023lowdim,he2026}, memorization and
generalization
\cite{achilli2025,biroli2024,latourellevigeant2026,li2024goodscore,lyu2025,shen2026,zhang2026wpo},
and replication in trained diffusion models
\cite{carlini2023,gu2025,somepalli2023}.

Compared with the existing literature discussed above, our main contributions
are the following two results on terminal localization and its critical scale
for empirical reverse diffusions with singular score errors.

First, we establish robustness of terminal localization under a
contraction-weighted error condition, including perturbations for which
full-path absolute continuity fails. Theorem~\ref{thm:main} proves almost-sure convergence to a training
point under progressively measurable, possibly singular score errors
satisfying the contraction-weighted condition. It also gives an explicit
trajectory bound around the same moving training point throughout a
terminal interval, with arbitrarily high probability.
Corollary~\ref{cor:positive-terminal-support} further shows that every
training point has positive terminal probability. When no deterministic
envelope is available, Proposition~\ref{prop:global-localization} and
Corollary~\ref{cor:mean-integrated} replace it by a realized weighted
integral along the perturbed path, respectively by its mean-integrated
version under the generated law, and still give almost-sure localization.
These conclusions allow
arbitrary initial probability laws and remain applicable when the exact and
perturbed full path laws are mutually singular (Lemma~\ref{lem:nonui}).

Second, we quantify the response to singular power-law errors and identify
the critical inverse-time scale at which terminal localization can fail.
Proposition~\ref{prop:matching} identifies the almost-sure leading response to the fixed-direction perturbation
\(e_r=M(T-r)^{-\gamma}\xi\), with \(M>0\), \(|\xi|=1\), and
\(1/2<\gamma<1\):
\[
 t^{\gamma-1}\bigl(X_t-y_{I_\infty}\bigr)
 \longrightarrow\frac{2M}{\gamma}\xi\qquad\text{a.s.}
\]
At the critical inverse-time scale, Proposition~\ref{prop:critical}
constructs perturbations of the form \(\kappa(T-r)^{-1}\xi\), activated on
a terminal interval, that produce a terminal limit outside the training set
with high probability even for arbitrarily small positive amplitudes, from
any initialization probability law.
These results quantify the robustness of empirical terminal attraction to
subcritical power-law errors and show how suitable critical perturbations can
alter the generated endpoint even at arbitrarily small amplitudes.

Next, we outline the main steps of the proof of terminal localization and
full terminal support. First, an exact
change of space and time gives the Gaussian
stochastic-localization representation
\cite{elalaoui2021localization,montanari2025sampling}, whose log-partition
function yields a global Lyapunov function. Its It\^o identity and the weighted
error condition force the distance to the training set to vanish almost
surely. Continuity and separation of the training points then give convergence
to one point, and hence high-probability entrance at sufficiently late
deterministic times. Second, once the trajectory is near a moving training
point, a local propagation estimate controls its displacement through the
geometric, score-error, and Brownian contributions. A first-exit argument
ensures that the trajectory stays near the same point and gives the explicit
bound throughout the remaining interval. Finally, truncated Girsanov
equivalence and an independent noise-control event show that every training
point has positive terminal probability, establishing full terminal support.

The remainder of the paper is organized as follows.
Section~\ref{sec:main} presents the full statements of the main results.
Section~\ref{sec:proof-localization} proves the quantitative terminal localization
theorem and the positive terminal support corollary.  Section~\ref{sec:rates} establishes the power-law asymptotics and the
critical perturbation result. Section~\ref{sec:conclusion} summarizes the contribution
and its scope. Appendix~\ref{app:prelim} establishes well-posedness before the
endpoint and terminal path-law singularity for an explicit infinite-energy
perturbation.

\section{Main results}\label{sec:main}

Under the model and assumption introduced in Section~\ref{sec:intro},
we have the following quantitative localization result.

\begin{theorem}[Quantitative localization and almost-sure convergence]\label{thm:main}
	Suppose that Assumption~\ref{ass:localization} holds.
	For every $\alpha\in(0,1)$, there exist a constant
	$b\in(0,T)$, an $\F_{T-b}$-measurable index $I_b\in\{1,\ldots,N\}$, and an event
	$\Hset_{b,\alpha}\in\F_{T-}$, with
	$\Prob(\Hset_{b,\alpha})\ge1-\alpha$, such that the following hold.
	\begin{enumerate}[label=\textup{(\roman*)},leftmargin=2.4em]
		\item \emph{Uniform terminal estimate.} On $\Hset_{b,\alpha}$, simultaneously
		for every $0<t\le b$,
		\begin{equation}\label{eq:main-envelope}
			\bigl|X_t-a_ty_{I_b}\bigr|\le \cR_{b,\alpha}(t)
			\triangleq G_b(t)+R_E(t,b)+\cN_{b,\alpha,d}(t).
		\end{equation}
		The three deterministic terms --- the geometric term $G_b$, the
		\emph{contraction-weighted response} $R_E$, and the Brownian envelope
		$\cN_{b,\alpha,d}$ --- are given by
		\begin{equation}\label{eq:Gb}
			G_b(t)\triangleq \sinh t\left[\frac{a_b\Lambda}{8\sinh b}+g(b)\right],
			\qquad
			g(b)\triangleq \frac{8(N-1)}{\Lambda}
			\exp\left[-\frac{\Lambda^2}{4(e^{2b}-1)}\right],
		\end{equation}
		and
		\begin{equation}\label{eq:error-response}
			\begin{aligned}
				R_E(t,b)&\triangleq 2\sinh t\int_t^b\frac{E(\sigma)}{\sinh \sigma}\dd \sigma ,\\
				\cN_{b,\alpha,d}(t)&\triangleq 4c(b)
				\sqrt{t\left[\ell_{\alpha,d}+2\log\bigl(1+\log_2(b/t)\bigr)\right]},\\
				c(b)&\triangleq \frac{\sinh b}{b},
				\qquad \ell_{\alpha,d}\triangleq \frac d2\log2+\log\frac{\pi^2}{3\alpha}.
			\end{aligned}
		\end{equation}
The function $\cR_{b,\alpha}$ is continuous on $(0,b]$ and tends to zero as $t\downarrow0$.
On $\Hset_{b,\alpha}$, $I_b$ is the unique minimizer of
$i\mapsto|X_t-a_ty_i|$ for every $0<t\le b$.
		\item \emph{Almost-sure convergence.} There exists an $\F_{T-}$-measurable
		index $I_\infty\in\{1,\ldots,N\}$ such that
		\[
		X_t\longrightarrow y_{I_\infty}\qquad\text{a.s. as }t\downarrow0.
		\]
		Moreover, $I_b=I_\infty$ almost surely on $\Hset_{b,\alpha}$.
		Almost surely, for all sufficiently small $t>0$, both
		$i\mapsto|X_t-a_ty_i|$ and $i\mapsto|X_t-y_i|$ have the unique
		minimizer $I_\infty$.
	\end{enumerate}
\end{theorem}

\begin{remark}[Interpretation of the terminal estimate]\label{rem:main-estimate}
The term $G_b$ bounds the initial displacement propagated from the start of
the terminal window and the contribution of the other mixture components.
The term $R_E$ bounds the accumulated effect of the score error, while
$\cN_{b,\alpha,d}$ bounds Brownian fluctuations throughout the window.

The terminal window is selected for the fixed
process using convergence of its entrance probability. The trajectory envelope
is explicit; we do not give a computable formula for choosing the window
uniformly over all errors with the same envelope.

Using the explicit expressions for these terms in Theorem~\ref{thm:main}, we obtain
the following convergence rates on \(\Hset_{b,\alpha}\) for power-law error
envelopes \(E(t)=O(t^{-\gamma})\), with \(0<\gamma<1\):
\[
\bigl|X_t-y_{I_b}\bigr|
=\begin{cases}
O\bigl(\sqrt{t\log\log(1/t)}\bigr),&0<\gamma\le1/2,\\
O(t^{1-\gamma}),&1/2<\gamma<1,
\end{cases}
\qquad(t\downarrow0).
\]
These \(O(\cdot)\) bounds hold for all \(0<t\le t_0\) on
\(\Hset_{b,\alpha}\), with the implicit constants and \(t_0\in(0,b]\)
chosen independently of the sample path.
\end{remark}

\begin{remark}[Denoiser-error criterion]\label{rem:denoiser}
The conditional mean \(\mu(x,t)\) of the moving centers, defined in
\eqref{eq:score}, is a denoiser: it estimates the
noiseless center from a noisy observation. The error condition can therefore
be restated in terms of this conditional mean.

Suppose Assumption~\ref{ass:localization}\textup{(\ref{ass:local})} holds. Write
\(D_E(t)\triangleq h_tE(t)\) for the deterministic envelope of the denoiser
correction \(\delta_r^\mu\triangleq h_{T-r}e_r\).
If
\[
 t\int_t^T\frac{D_E(\sigma)}{\sigma^2}\dd \sigma \longrightarrow0
 \qquad(t\downarrow0),
\]
then all conclusions of Theorem~\ref{thm:main} hold.
In particular, this criterion holds if \(D_E(t)\to0\) as \(t\downarrow0\).

This reformulation is the natural one for learned scores: the denoiser correction $h_{T-r}e_r$ is the object that appears in the drift, and a bound $\|s_\theta-s\|_{L^2(w(t))}=O(t^{-\gamma})$ in $L^2(w(t))$ translates into the bound $O(t^{1-\gamma})$ for the correction $D_E$. It is a pathwise sufficient criterion, not a training-risk guarantee; see \cite[\S4]{liu2025} for the relative-entropy viewpoint.
\end{remark}
The following example shows that the weighted condition is strictly weaker
than the integrated squared-error condition
\[
t\int_t^T E(\sigma)^2\dd\sigma\longrightarrow0
\qquad(t\downarrow0),
\]
and permits sufficiently narrow supercritical spikes.
\begin{example}[Sparse supercritical spikes]\label{rem:spikes}
Cauchy--Schwarz gives
\[
 t\int_t^T\frac{E(\sigma)}\sigma \dd \sigma 
 \le\left(t\int_t^T E(\sigma)^2\dd \sigma \right)^{1/2}
       \left(1-\frac tT\right)^{1/2}
 \le\left(t\int_t^T E(\sigma)^2\dd \sigma \right)^{1/2}.
\]
Thus the squared-error condition implies \eqref{eq:weighted-vanish}, but
the converse fails. Take $T=1$, $\sigma_k=2^{-k}$ for $k\ge2$, and
\[
 A=\bigcup_{k\ge2}[\sigma_k,\sigma_k+\sigma_k^3],
 \qquad E(t)=t^{-2}\one_A(t),\quad0<t\le1.
\]
The intervals are disjoint and only finitely many meet $[\delta,1]$ for any
$\delta>0$, so $E\in L^2([\delta,1])$. Each spike contributes at most one
to the weighted integral:
\[
 \int_{\sigma_k}^{\sigma_k+\sigma_k^3}\frac{E(\sigma)}\sigma \dd \sigma 
 =\int_{\sigma_k}^{\sigma_k+\sigma_k^3}\sigma^{-3}\dd \sigma \le1.
\]
There are $O(1+\log(1/t))$ full spikes above $t$ and at most one partial
spike crossing $t$. Hence
\[
 t\int_t^1\frac{E(\sigma)}\sigma \dd \sigma 
 =O\bigl(t[1+\log(1/t)]\bigr)\longrightarrow0\qquad(t\downarrow0).
\]
On the other hand,
\[
\begin{aligned}
 \sigma_k\int_{\sigma_k}^{\sigma_k+\sigma_k^3}E(\sigma)^2\dd \sigma
 &=\frac{\sigma_k}{3}\left[\frac{1}{\sigma_k^3}
 -\frac{1}{(\sigma_k+\sigma_k^3)^3}\right]\\
 &=\frac{3+3\sigma_k^2+\sigma_k^4}{3(1+\sigma_k^2)^3}
 \longrightarrow1\qquad(k\to\infty).
\end{aligned}
\]
Since $\sigma_k\int_{\sigma_k}^1 E(\sigma)^2\dd \sigma $ is at least this contribution, the
squared-error condition fails. The deterministic perturbation
$e_r=E(1-r)\xi$, for a fixed unit vector $\xi$, therefore satisfies the
localization criterion despite these supercritical spikes and infinite
squared-error energy. The inverse-time threshold concerns persistent
power-law perturbations, not a pointwise upper bound on the error.
\end{example}

Almost-sure convergence gives a terminal law supported on $S$.
The next result shows that it assigns positive mass to every target atom.

\begin{corollary}[Positive terminal selection probabilities and full support]
    \label{cor:positive-terminal-support}
    Under the assumptions of Theorem~\ref{thm:main}, let $I_\infty$ be its terminal
    index. Define $X_0\triangleq\lim_{t\downarrow0}X_t=y_{I_\infty}$ as its
    almost-sure continuous extension to the generation endpoint $t=0$, and put
    $v_0^*\triangleq\cL(X_0)$, the limiting generated law. Then
    \begin{equation}\label{eq:terminal-law}
        \begin{gathered}
            q_i\triangleq\Prob(I_\infty=i)>0\quad(1\le i\le N),
            \qquad \sum_{i=1}^Nq_i=1,\\
            v_0^*=\sum_{i=1}^Nq_i\delta_{y_i},
            \qquad \operatorname{supp}(v_0^*)=S.
        \end{gathered}
    \end{equation}
    We call $q_i$ the terminal selection probability of $y_i$ (equivalently, its
    terminal sampling weight). Support recovery refers only to the identity
    $\operatorname{supp}(v_0^*)=\operatorname{supp}(u_0)=S$; it is a geometric
    statement and by itself does not control the $q_i$.
\end{corollary}
Its proof is given in Section~\ref{sec:positive-terminal-support}.

\begin{remark}[Localization and distributional fidelity]\label{rem:fidelity}
Terminal localization specifies where a trajectory ends. Distributional
fidelity additionally requires the terminal selection probabilities to
approximate their empirical values. With the convention
$\|\nu-\mu\|_{\mathrm{TV}}=\sup_A|\nu(A)-\mu(A)|$, the terminal law in
\eqref{eq:terminal-law} satisfies
\begin{equation}\label{eq:terminal-tv}
 \|v_0^*-u_0\|_{\mathrm{TV}}
 =\frac12\sum_{i=1}^N\left|q_i-\frac1N\right|.
\end{equation}
Almost-sure convergence to $S$ and positivity of the $q_i$ do not supply
a small bound on this quantity. Theorem~\ref{thm:main} controls spatial
displacement near the selected atom; it does not estimate selection-weight
errors. The distinction itself is not new --- it is already drawn in \cite[\S1]{liu2025} and
implicitly in the memorization literature --- but it is worth isolating because
localization is often read as a statement about the sampling distribution. Example~\ref{ex:changed-weights} below exhibits this distinction at
the same initialization as the exact sampler.
\end{remark}

\begin{example}[Preserved support with changed selection probabilities]
\label{ex:changed-weights}
Let $d=1$, $S=\{-1,1\}$, and fix $0<b<T$ and $p\in(0,1)$.
Write $\theta=2p-1$ and
\[
 f_{\pm,t}(x)=\varphi_{h_t}(x\mp a_t),\qquad
 w_p(x,t)=p f_{+,t}(x)+(1-p)f_{-,t}(x),\qquad
 m_t(x)=\tanh(a_tx/h_t).
\]
Thus $w_{1/2}=w$. Retain the exact sampler's initialization
$X_T\sim w(T)$ and use the error $e_r=e^{p,b}(\widehat X_r,T-r)$, where
\begin{equation}\label{eq:two-atom-error}
 e^{p,b}(x,t)=\one_{[b,T]}(t)\frac{a_t}{h_t}
 \left[\tanh\!\left(\frac{a_tx}{h_t}+\frac12\log\frac p{1-p}\right)
             -\tanh\!\left(\frac{a_tx}{h_t}\right)\right].
\end{equation}
Generation uses the score $\partial_x\log w_p$ down to $t=b$, then
restores the original score for $0<t<b$. This error satisfies
Assumption~\ref{ass:localization}, with the bounded all-space envelope
\begin{equation}\label{eq:two-atom-envelope}
 E_{p,b}(t)=\one_{[b,T]}(t)\frac{\kappa_p}{\sinh t},\qquad
 \kappa_p=\tanh\!\left(\frac14\left|\log\frac p{1-p}\right|\right).
\end{equation}
Its terminal law is
$q_+(p,b)\delta_1+(1-q_+(p,b))\delta_{-1}$, with
$0<q_+(p,b)<1$. The function $q_+(p,b)$ is strictly increasing in $p$,
and $q_+(1/2,b)=1/2$. Hence $p>1/2$ gives a different terminal law with
exactly the same two-point support and the same initialization as the
exact sampler.

For $T=1$, $b=1/2$ and $p=3/4$, one has
$\kappa_p=2-\sqrt3$ and
$\sup_t E_{p,b}(t)=(2-\sqrt3)/\sinh(1/2)\simeq0.51420$.
The exact integral formula
\eqref{eq:two-atom-weight} below gives numerically
\[
 q_+(3/4,1/2)\simeq0.56498,\qquad
 \|v_0^*-u_0\|_{\mathrm{TV}}\simeq0.06498.
\]
These decimals are illustrative: composite Simpson quadrature on two
standardized Gaussian variables in $[-9,9]$, with $200$ and $400$
subintervals per variable, agrees at the displayed precision. Strict
inequality and support preservation follow from the proof, independently
of the numerical calculation.
\end{example}
Its proof and explicit Gaussian integral are given in
Section~\ref{sec:two-atom-proof}.

\begin{remark}[Path laws and terminal position laws]\label{rem:path-laws}
For the deterministic perturbation $e_r=M(T-r)^{-\gamma}\xi$, with
$M>0$, $|\xi|=1$ and $1/2\le\gamma<1$, let $P$ and $Q$ be the exact and
perturbed laws on the common canonical path space of Appendix~\ref{app:prelim},
both started from $w(T)$. Writing $P^{\mathrm{end}}$ and $Q^{\mathrm{end}}$
for their terminal position laws and letting $\F_{T-}^0$ denote the
sigma-field generated by the increasing-time trajectory $\widehat X_r$ for $r<T$,
Lemma~\ref{lem:nonui} and Corollary~\ref{cor:positive-terminal-support} give
\[
\begin{gathered}
P\perp Q\quad\text{on }\F_{T-}^0,\qquad
P^{\mathrm{end}}\sim Q^{\mathrm{end}},\\
P^{\mathrm{end}}=u_0,\qquad
Q^{\mathrm{end}}=\sum_{i=1}^Nq_i\delta_{y_i},\quad q_i>0.
\end{gathered}
\]
Thus singular full path laws can have equivalent terminal position laws.
\end{remark}

\begin{remark}[Scope of the conclusions]\label{rem:scope}
The deterministic error bound is imposed along the perturbed trajectory.
Proposition~\ref{prop:global-localization} below replaces it, for qualitative
almost-sure localization, by a realized weighted-integral condition, and
Corollary~\ref{cor:mean-integrated} gives a sufficient condition averaged
under the generated law $v(t)$. An $L^2(w(t))$ estimate under the reference
law does not by itself imply any of these trajectory or generated-law
conditions. The sample index is chosen using the information
available at the start of the terminal window; it need not agree with the
terminal index of a coupled exact trajectory. The point set is fixed:
constants and admissible windows depend on its geometry, including $\Lambda$.
We do not claim estimates uniform in the number of points or as their
separation tends to zero.
The localization, full-support, subcritical power-law, and critical displacement
results allow arbitrary initial probability laws, without entropy or moment
bounds.

Concretely, with $g$ as in \eqref{eq:Gb}, the first requirement in
\eqref{eq:conditions}, $\sinh(b)g(b)\le a_b\Lambda/32$, forces
\[
 b\lesssim\frac{\Lambda^2}{8\log\bigl(256(N-1)/\Lambda^2\bigr)},
\]
so the admissible window is shorter for larger $N$ or smaller $\Lambda$.
The remaining requirements, in particular $\varphi_b\le\eta$, are also
satisfied for all sufficiently small $b$ (Lemma~\ref{lem:energy}); no lower
bound on $b$ is needed. The dependence on $N$ and $\Lambda$ thus enters only
through the upper bound above.

\end{remark}

For the following fixed-direction perturbation, the integral response determines
an explicit leading displacement.
\begin{proposition}[Leading-order response to a power-law perturbation]\label{prop:matching}
	Fix $M>0$,
	$1/2<\gamma<1$ and a unit vector $\xi\in\R^d$, and let the error be
	\begin{equation}\label{eq:power-error}
		e_r\triangleq M(T-r)^{-\gamma}\xi,\qquad 0\le r<T .
	\end{equation}
	Let $I_\infty$ be the terminal index of Theorem~\ref{thm:main}. Then
	\begin{equation}\label{eq:matching}
		t^{\gamma-1}\bigl(X_t-y_{I_\infty}\bigr)\ \longrightarrow\ \frac{2M}{\gamma}\,\xi
	\end{equation}
	almost surely as $t\downarrow0$.
\end{proposition}

The restriction $\gamma>1/2$ in Proposition~\ref{prop:matching} makes the
Brownian contribution lower order than the deterministic error response.
Terminal localization still holds for $\gamma\leq1/2$, where that dominance
is not asserted.

At the critical scale, we turn on the perturbation only during a terminal
interval. Before this interval, the process evolves with the exact score.
\begin{proposition}[Critical terminal displacement]\label{prop:critical}
	Let $v_T$ be any initialization probability measure on $\R^d$.
	Fix a unit vector
	$\xi\in\R^d$, $0<\kappa<\Lambda/64$ and $0<\alpha<1$, and let
	\begin{equation}\label{eq:critical-error}
		e_r^{\kappa,b}\triangleq \frac{\kappa}{T-r}\,\xi\,\one_{[T-b,T)}(r),
	\end{equation}
	where $b\in(0,T)$ and the switch time $T-b$ is part of the definition of the error. For each $b$,
	let $\widehat X^{\kappa,b}$ denote the solution of \eqref{eq:main-sde} with
	error $e^{\kappa,b}_r$, using the same initial random variable
	$\widehat X_0\sim v_T$ and Brownian motion $B$, and set
	$X_t^{\kappa,b}\triangleq\widehat X_{T-t}^{\kappa,b}$. Then there exist a deterministic
	$b=b(v_T,\kappa,\alpha)\in(0,T)$, with the model parameters fixed,
	an $\F_{T-b}$-measurable index $I_b\in\{1,\ldots,N\}$, and an event
	$\Hset^{\mathrm{crit}}_{b,\alpha}\in\F_{T-}$ with
	$\Prob(\Hset^{\mathrm{crit}}_{b,\alpha})\ge1-\alpha$ such that, on this event,
	simultaneously for every $0<t\le b$,
	\[
		\bigl|X^{\kappa,b}_t-a_ty_{I_b}\bigr|\le\frac{a_t\Lambda}{4}.
	\]
	Moreover, for $\Prob$-almost every path of $\Hset^{\mathrm{crit}}_{b,\alpha}$,
	\begin{equation}\label{eq:critical-limit}
		\begin{gathered}
			X^{\kappa,b}_t\longrightarrow y_{I_b}+2\kappa\xi
			\qquad(t\downarrow0),\\
			\dist\bigl(y_{I_b}+2\kappa\xi,S\bigr)=2\kappa>0.
		\end{gathered}
	\end{equation}
\end{proposition}
The perturbation \eqref{eq:critical-error} is switched on only on $[T-b,T)$ so that,
before $T-b$, the perturbed solution coincides pathwise with the exact one; the
almost-sure convergence of that exact process then supplies the entrance estimate for
every initialization law. Inside the window the accumulated response tends to $2\kappa$
as $t\downarrow0$, whatever the size of $\kappa$, and \eqref{eq:critical-weighted}
shows that \eqref{eq:weighted-vanish} fails for these perturbations.
\begin{remark}[Dependence on the switching time]\label{rem:rates}
    Since $c(b),a_b\to1$ as $b\downarrow0$, the restriction
    $\kappa<\Lambda/64$ ensures $2\kappa c(b)\le a_b\Lambda/32$
    for all sufficiently small $b$, as required in \eqref{eq:crit-conditions}.
	In Proposition~\ref{prop:critical}, changing $b$ changes both the perturbation
	and the process, so the high-probability conclusions for varying $b$ do not
	imply an almost-sure statement for a fixed process. Nevertheless, all these
	processes coincide before their respective switch times with one fixed
	unperturbed process. Its almost-sure convergence supplies the entrance
	probability needed to choose $b$, for any $v_T$.
	Since $\kappa$ can be
	arbitrarily small, the proposition excludes a universal constant $\kappa_*>0$
	guaranteeing almost-sure convergence to $S$ for every perturbation satisfying
	$|e_r|\le\kappa_*/(T-r)$. It neither classifies all critical perturbations nor asserts that the weighted
	condition \eqref{eq:weighted-vanish} is necessary for convergence to $S$;  compare
	Example~\ref{rem:spikes} for supercritical errors that do localize.
\end{remark}

\section{Proofs of localization and positivity}\label{sec:proof-localization}

We first establish almost-sure convergence through a global Lyapunov
identity after an exact change of space and time. This gives entrance near a
moving center with probability tending to one at deterministic times.
The local propagation, error and noise estimates then rule out a first exit
and provide an explicit bound throughout a deterministic terminal window.
We finally show that each center is selected with positive probability.

\subsection{Global convergence and deterministic-time entrance}\label{sec:entrance}

For deterministic $0<b<T$, define
\begin{equation}\label{eq:entry}
    \rho_b^{\mathrm{ent}}\triangleq\frac{a_b\Lambda}{8},\qquad
    I_b\triangleq\min\argmin_{1\le i\le N}|\widehat X_{T-b}-a_by_i|.
\end{equation}
The smallest index is selected on a tie. The entrance event is
\begin{equation}\label{eq:entrance-event}
    \Eset_b\triangleq\{\min_i|\widehat X_{T-b}-a_by_i|\le \rho_b^{\mathrm{ent}}\}\in\F_{T-b}.
\end{equation}
The index $I_b$ is $\F_{T-b}$-measurable and unique on $\Eset_b$, since
the entrance balls have center separation $a_b\Lambda=8\rho_b^{\mathrm{ent}}$.
Write \(\varphi_t\triangleq\Prob(\Eset_t^{\mathrm c})\) for \(0<t<T\), the
entrance failure probability at time \(t\); in particular
\(\varphi_b=\Prob(\Eset_b^{\mathrm c})\).
That this index remains nearest throughout the terminal window will follow
from the first-exit argument.

\begin{proposition}[Localization under a realized weighted-integral condition]
\label{prop:global-localization}
In the model~\eqref{eq:main-sde}, assume that $e$ is progressively
measurable, that the compatible continuous adapted solutions exist on every
$[0,T-\delta]$, $0<\delta<T$, and that
\[
 \int_0^{T-\delta}|e_r|^2\dd r<\infty\quad\text{a.s.}
 \qquad(0<\delta<T).
\]
A deterministic error envelope is not required. Suppose instead that
\begin{equation}\label{eq:realized-weighted}
 A_e(t)\triangleq t\int_t^T\frac{|e_{T-\sigma}|}{\sigma}\dd\sigma
 \longrightarrow0\qquad\text{a.s. as }t\downarrow0.
\end{equation}
Then there exists an $\F_{T-}$-measurable index
$I_\infty\in\{1,\ldots,N\}$ such that
\[
 X_t\longrightarrow y_{I_\infty}
 \qquad\text{a.s. as }t\downarrow0.
\]
This conclusion requires no initial entropy or moment bound. In particular,
it applies under Assumption~\ref{ass:localization}.
\end{proposition}

The space--time transformation below identifies the exact reverse
Ornstein--Uhlenbeck dynamics with a Gaussian stochastic-localization process
for the empirical measure, up to deterministic rescaling
\cite{elalaoui2021localization,montanari2025sampling}.
We use this representation to study progressively measurable perturbations
whose accumulated magnitude in the transformed time is sublinear.
The proof combines the associated log-partition function with a localized
martingale argument to establish almost-sure convergence to a training point.

\begin{proof}[Proof of Proposition~\ref{prop:global-localization}]
\stepheading{Space--time change and a Lyapunov identity.}
Set
\begin{equation}\label{eq:global-time-change}
\begin{gathered}
q(t)\triangleq\coth t-1=\frac{a_t}{\sinh t},\qquad
q_0\triangleq q(T)>0,\\
t(q)\triangleq\frac12\log\left(1+\frac2q\right),\qquad
r(q)\triangleq T-t(q),\qquad q\ge q_0,\\
Z_q\triangleq\frac{\widehat X_{r(q)}}{\sinh t(q)}.
\end{gathered}
\end{equation}
The time-changed filtration is
$\mathcal G_q\triangleq\F_{r(q)}$.
Since $\mathrm d q/\mathrm d r=\sinh^{-2}(T-r)$, the process
\[
\widetilde B_q\triangleq
\int_0^{r(q)}\frac{\dd B_s}{\sinh(T-s)},\qquad q\ge q_0,
\]
has covariations
\[
\langle\widetilde B^{(j)},\widetilde B^{(k)}\rangle_q
=\delta_{jk}\int_0^{r(q)}\frac{\dd s}{\sinh^2(T-s)}
=\delta_{jk}(q-q_0)
\]
and starts at zero. L\'evy's characterization therefore shows that
$(\widetilde B_{q_0+\lambda})_{\lambda\ge0}$ is a standard
$d$-dimensional Brownian motion in the filtration
$(\mathcal G_{q_0+\lambda})_{\lambda\ge0}$.

The identities $h_t=2a_t\sinh t$ and
$1-2/h_t=-\coth t$ rewrite the original drift as
\[
\widehat X_r+2s(\widehat X_r,t)
=-\coth t\,\widehat X_r+\frac1{\sinh t}
  \sum_{i=1}^N \omega_i(\widehat X_r,t)y_i,\qquad t=T-r.
\]
The spatial rescaling has derivative
\[
\frac{\dd}{\dd r}\frac1{\sinh(T-r)}
=\frac{\cosh(T-r)}{\sinh^2(T-r)}
=\frac{\coth t}{\sinh t},\qquad t=T-r.
\]
Substituting the drift above into \eqref{eq:main-sde} gives
\[
\begin{aligned}
\dd\left(\frac{\widehat X_r}{\sinh t}\right)
&=\frac1{\sinh t}\dd\widehat X_r
  +\frac{\coth t}{\sinh t}\widehat X_r\dd r\\
&=\left[\frac1{\sinh^2 t}\sum_{i=1}^N\omega_i(\widehat X_r,t)y_i
  +\frac{2e_r}{\sinh t}\right]\dd r
  +\frac{\sqrt2}{\sinh t}\dd B_r.
\end{aligned}
\]
For $x=\sinh t\,z$ and $q=a_t/\sinh t$, the identity
$h_t=2a_t\sinh t$ gives
\[
-\frac{|x-a_ty_i|^2}{2h_t}
=-\frac{|x|^2}{2h_t}+\frac{x\cdot y_i}{2\sinh t}
  -\frac{a_t|y_i|^2}{4\sinh t}
=-\frac{|x|^2}{2h_t}+\frac{z\cdot y_i}{2}-\frac{q|y_i|^2}{4}.
\]
Hence
\[
\begin{aligned}
\omega_i(x,t)
&=\frac{\exp(-|x-a_ty_i|^2/(2h_t))}
{\sum_{j=1}^N\exp(-|x-a_ty_j|^2/(2h_t))}\\
&=\frac{\exp(z\cdot y_i/2-q|y_i|^2/4)}
{\sum_{j=1}^N\exp(z\cdot y_j/2-q|y_j|^2/4)}
\triangleq\pi_i(q,z).
\end{aligned}
\]
Using
\[
\frac{\dd r(\zeta)}{\dd\zeta}=\sinh^2 t(\zeta),\qquad
\omega_i(\widehat X_{r(\zeta)},t(\zeta))=\pi_i(\zeta,Z_\zeta),
\]
we obtain
\[
\begin{aligned}
Z_q-Z_{q_0}
&=\int_0^{r(q)}\left[
\frac{\sum_{i=1}^N\omega_i(\widehat X_s,T-s)y_i}{\sinh^2(T-s)}
+\frac{2e_s}{\sinh(T-s)}\right]\dd s\\
&\quad+\sqrt2\int_0^{r(q)}\frac{\dd B_s}{\sinh(T-s)}\\
&=\int_{q_0}^q\left[
\sum_{i=1}^N\pi_i(\zeta,Z_\zeta)y_i
+2\sinh t(\zeta)e_{r(\zeta)}\right]\dd\zeta
+\sqrt2\,\widetilde B_q.
\end{aligned}
\]
In differential form,
\begin{equation}\label{eq:global-transformed-sde}
\dd Z_q=\left[\sum_{i=1}^N\pi_i(q,Z_q)y_i
+2\sinh t(q)\,e_{r(q)}\right]\dd q
+\sqrt2\,\dd\widetilde B_q.
\end{equation}
For $q>0$ and $z\in\R^d$, define
\begin{equation}\label{eq:global-lyapunov}
\begin{aligned}
\mathcal V(q,z)
&\triangleq\log\left[\frac1N\sum_{i=1}^N
\exp\left(-\frac q4\left|y_i-\frac zq\right|^2\right)\right]\\
&=\log\left[\frac1N\sum_{i=1}^N
e^{z\cdot y_i/2-q|y_i|^2/4}\right]-\frac{|z|^2}{4q}
\le0.
\end{aligned}
\end{equation}
Set
\[
\operatorname{Var}_\pi(y)\triangleq
\sum_{i=1}^N\pi_i(q,z)|y_i|^2
-\left|\sum_{i=1}^N\pi_i(q,z)y_i\right|^2.
\]
Direct differentiation gives
\[
\begin{aligned}
\partial_q\mathcal V
&=-\frac14\sum_{i=1}^N\pi_i(q,z)|y_i|^2+\frac{|z|^2}{4q^2},
&\Delta_z\mathcal V
&=\frac14\operatorname{Var}_\pi(y)-\frac d{2q},\\
\nabla_z\mathcal V
&=\frac12\left(\sum_{i=1}^N\pi_i(q,z)y_i-\frac zq\right).
\end{aligned}
\]
Applying It\^o's formula to \eqref{eq:global-transformed-sde}, first
on bounded $q$-intervals and then by localization, gives the exact identity
\begin{equation}\label{eq:global-ito}
\begin{aligned}
\dd\mathcal V(q,Z_q)
&=\left[\frac14\left|\sum_{i=1}^N\pi_i(q,Z_q)y_i-\frac{Z_q}{q}\right|^2
-\frac d{2q}\right]\dd q\\
&\quad+\left(\sum_{i=1}^N\pi_i(q,Z_q)y_i-\frac{Z_q}{q}\right)
\cdot\left[e_{r(q)}\,\sinh t(q)\dd q
+\frac1{\sqrt2}\dd\widetilde B_q\right].
\end{aligned}
\end{equation}
\proofstep{Sublinear growth of the martingale.}
Let
\[
M_q\triangleq\frac1{\sqrt2}\int_{q_0}^q
\left(\sum_{i=1}^N\pi_i(\zeta,Z_\zeta)y_i-\frac{Z_\zeta}{\zeta}\right)
\cdot\dd\widetilde B_\zeta.
\]
We claim that
\[
\frac{M_q}{q}\longrightarrow0\qquad\text{a.s. as }q\to\infty.
\]

To prove this, integrate \eqref{eq:global-transformed-sde} and divide by $q$:
\[
\frac{Z_q}{q}=\frac{Z_{q_0}}q
+\frac1q\int_{q_0}^q\sum_{i=1}^N\pi_i(\zeta,Z_\zeta)y_i\dd\zeta
+\frac2q\int_{q_0}^q\sinh t(\zeta)e_{r(\zeta)}\dd\zeta
+\frac{\sqrt2\,\widetilde B_q}q.
\]
Since $Z_{q_0}$ is finite almost surely and
$|\sum_i\pi_i(q,z)y_i|\le R$, it remains to control the forcing term.
Work on a single probability-one event on which the local error integrals
are finite and \eqref{eq:realized-weighted} holds. Such an event is obtained
by a countable exhaustion of the truncated intervals. Put $t=t(q)\downarrow0$
and $c(T)=\sinh T/T$. Changing variables to $\sigma=t(\zeta)$ yields
\begin{equation}\label{eq:global-forcing}
\begin{aligned}
\frac2q\int_{q_0}^q\sinh t(\zeta)\,|e_{r(\zeta)}|\dd \zeta
&=\frac{2\sinh t}{a_t}\int_t^T\frac{|e_{T-\sigma}|}{\sinh\sigma}\dd\sigma\\
&\le2e^t c(T)t\int_t^T\frac{|e_{T-\sigma}|}{\sigma}\dd \sigma
=2e^t c(T)A_e(t)\longrightarrow0,
\qquad(q\to\infty).
\end{aligned}
\end{equation}
The last line uses $\sinh\sigma\ge\sigma$, $\sinh t\le c(T)t$, and
\eqref{eq:realized-weighted}.

The Brownian strong law and the integrated equation therefore imply
\[
\limsup_{q\to\infty}\frac{|Z_q|}{q}\le R\qquad\text{a.s.}
\]
By continuity on bounded $q$-intervals,
\begin{equation}\label{eq:global-boundedness}
\sup_{q\ge q_0}\frac{|Z_q|}{q}<\infty\qquad\text{a.s.}
\end{equation}

For each positive integer $K$, set
\[
\tau_K\triangleq\inf\{q\ge q_0:\frac{|Z_q|}{q}\ge K\},
\qquad\inf\varnothing=\infty.
\]
Before $\tau_K$ the integrand of $M$ has norm at most $(R+K)/\sqrt2$,
so It\^o's isometry gives
\[
\E|M_{q\wedge\tau_K}|^2
\le\frac{(R+K)^2}{2}(q-q_0).
\]
For $\varepsilon>0$ and integers $j$ with $2^j\ge q_0$,
Markov's and Doob's inequalities yield
\[
\begin{aligned}
\Prob\left(\sup_{q_0\le\zeta\le2^{j+1}}|M_{\zeta\wedge\tau_K}|
>\varepsilon2^j\right)
&\le\frac{4\E|M_{2^{j+1}\wedge\tau_K}|^2}{\varepsilon^2 2^{2j}}\\
&\le\frac{4(R+K)^2}{\varepsilon^2}2^{-j}.
\end{aligned}
\]
These probabilities are summable. Borel--Cantelli, applied simultaneously
to all integers $K\ge1$ and $\varepsilon\in\mathbb Q_{>0}$, therefore gives
\[
\sup_{2^j\le q\le2^{j+1}}\frac{|M_{q\wedge\tau_K}|}{q}\longrightarrow0
\quad\text{a.s. for every }K\ge1.
\]
Intersecting with the probability-one event in \eqref{eq:global-boundedness},
choose an integer $K>\sup_{q\ge q_0}|Z_q|/q$. Then $\tau_K=\infty$,
proving $M_q/q\to0$ almost surely.

\proofstep{Terminal localization.}
Integrating \eqref{eq:global-ito} gives
\[
\begin{aligned}
\mathcal V(q,Z_q)
&=\mathcal V(q_0,Z_{q_0})
+\frac14\int_{q_0}^q\left|\sum_{i=1}^N\pi_i(\zeta,Z_\zeta)y_i
-\frac{Z_\zeta}{\zeta}\right|^2\dd\zeta-\frac d2\log(q/q_0)\\
&\quad+\int_{q_0}^q\left(\sum_{i=1}^N\pi_i(\zeta,Z_\zeta)y_i
-\frac{Z_\zeta}{\zeta}\right)\cdot e_{r(\zeta)}\,\sinh t(\zeta)\dd\zeta+M_q.
\end{aligned}
\]
On the pathwise boundedness event,
\[
\begin{aligned}
&\left|\int_{q_0}^q\left(\sum_{i=1}^N\pi_i(\zeta,Z_\zeta)y_i-\frac{Z_\zeta}{\zeta}\right)
\cdot e_{r(\zeta)}\,\sinh t(\zeta)\dd \zeta\right|\\
&\qquad\le\left(R+\sup_{\zeta\ge q_0}\frac{|Z_\zeta|}{\zeta}\right)
\int_{q_0}^q\sinh t(\zeta)\,|e_{r(\zeta)}|\dd \zeta\\
&\qquad=o(q)
\end{aligned}
\]
by \eqref{eq:global-forcing}.
Using $M_q=o(q)$ and discarding the nonnegative quadratic term gives
\[
\mathcal V(q,Z_q)\ge
\mathcal V(q_0,Z_{q_0})-\frac d2\log(q/q_0)+o(q).
\]
Since $\mathcal V\le0$,
\[
0\ge\frac{\mathcal V(q,Z_q)}q
\ge\frac{\mathcal V(q_0,Z_{q_0})}q
-\frac{d}{2q}\log(q/q_0)+o(1)\longrightarrow0
\qquad\text{a.s.}
\]
Every summand in \eqref{eq:global-lyapunov} is at most
$\exp[-q\dist(Z_q/q,S)^2/4]$; hence
\begin{equation}\label{eq:global-distance}
0\le\dist\left(\frac{Z_q}{q},S\right)^2
\le-\frac4q\mathcal V(q,Z_q)\longrightarrow0
\qquad\text{a.s.}
\end{equation}

On the full-probability event in \eqref{eq:global-distance}, the continuous
path $q\mapsto Z_q/q$ eventually lies in
$\bigcup_{i=1}^NB(y_i,\Lambda/3)$.
These balls are pairwise disjoint.
The image of the connected interval $[q_*,\infty)$ under a continuous
path is connected, so once the path remains in this union it belongs
to a single one of the balls.
Within that ball the corresponding $y_i$ is the nearest point of $S$;
\eqref{eq:global-distance} therefore implies $\frac{Z_q}{q}\to y_i$.
Since $X_t=\sinh t\,Z_{q(t)}=a_tZ_{q(t)}/q(t)$ and $a_t\to1$,
also $X_t\to y_i$ as $t\downarrow0$.

To specify a measurable terminal index, fix deterministic rational times
$r_n\in[0,T)$ with $r_n\uparrow T$ and set
\[
J_n\triangleq\min\argmin_{1\le i\le N}|\widehat X_{r_n}-y_i|,
\qquad I_\infty\triangleq\liminf_{n\to\infty}J_n.
\]
Each $J_n$ is $\F_{r_n}$-measurable, so $I_\infty\in\{1,\ldots,N\}$
is $\F_{T-}$-measurable. On the full-probability convergence event above,
separation of $S$ makes $J_n$ eventually equal to the index of the limit.
Thus $X_t\to y_{I_\infty}$ almost surely. This definition uses only countably
many coordinate evaluations and does not require completion of the path
sigma-field.
\end{proof}

\begin{corollary}[A mean-integrated sufficient condition]
\label{cor:mean-integrated}
Retain the local solution and integrability assumptions of
Proposition~\ref{prop:global-localization}. For some $0<b<T$, set
$t_n=b2^{-n}$, $n\ge0$. If
\begin{equation}\label{eq:mean-integrated}
 \sum_{n=0}^\infty t_n\int_{t_n}^b
 \frac{\E|e_{T-\sigma}|}{\sigma}\dd\sigma<\infty,
\end{equation}
then $X_t$ converges almost surely to a training point. In particular,
this holds if, for some $C<\infty$ and $0\le\gamma<1$,
\begin{equation}\label{eq:mean-power}
 \E|e_{T-t}|\le Ct^{-\gamma}\qquad\text{for a.e. }0<t\le b.
\end{equation}
All expectations are under the actual perturbed process. For a Markov
score error, the integrand in~\eqref{eq:mean-integrated} is
\[
 \E|e_{T-t}|=\int_{\R^d}|s_\theta(x,t)-s(x,t)|\,v(t)(\dd x).
\]
By Cauchy--Schwarz, an $L^2(v(t))$ score-error bound of the form
$\bigl(\int|s_\theta-s|^2\,\dd v(t)\bigr)^{1/2}\le Ct^{-\gamma}$
also implies~\eqref{eq:mean-power}.
\end{corollary}
\begin{proof}[Proof of Corollary \ref{cor:mean-integrated}]
Put $J_n=t_n\int_{t_n}^b|e_{T-\sigma}|/\sigma\dd\sigma$.
Tonelli's theorem and~\eqref{eq:mean-integrated} give
$\E\sum_nJ_n<\infty$, so $J_n\to0$ almost surely. For
$t_{n+1}\le t\le t_n$, nonnegativity gives
\[
 t\int_t^b\frac{|e_{T-\sigma}|}{\sigma}\dd\sigma
 \le t_n\int_{t_{n+1}}^b\frac{|e_{T-\sigma}|}{\sigma}\dd\sigma
 =2J_{n+1}.
\]
The integral from $b$ to $T$ is finite almost surely by local square
integrability, so its product with $t$ tends to zero. Thus
\eqref{eq:realized-weighted} holds and
Proposition~\ref{prop:global-localization} applies. Under
\eqref{eq:mean-power}, for $\gamma>0$ one has
$\E J_n\le(C/\gamma)t_n^{1-\gamma}$; for $\gamma=0$ one has
$\E J_n\le Ct_n\log(b/t_n)$. Both bounds are summable.
\end{proof}

\begin{remark}[Scope of the integrated extension]\label{rem:integrated-scope}
Proposition~\ref{prop:global-localization} and
Corollary~\ref{cor:mean-integrated} give qualitative almost-sure
localization. For example, let $Z$ be an integrable $\F_0$-measurable random variable
that is not essentially bounded and let $\xi$ be a fixed unit vector.
The error $e_r=Z(T-r)^{-\gamma}\xi$, $0\le\gamma<1$, satisfies
\eqref{eq:mean-power} but has no finite deterministic envelope of the form
\eqref{eq:envelope}. The explicit deterministic envelope in Theorem~\ref{thm:main}
and the positive mass of every atom in
Corollary~\ref{cor:positive-terminal-support} retain
Assumption~\ref{ass:localization}. In particular, a reference-law learning
bound under $w(t)$ cannot be substituted for the generated-law condition
\eqref{eq:mean-integrated} without an additional comparison argument.
\end{remark}

\begin{lemma}[Vanishing of the entrance failure probability]\label{lem:energy}
Under Assumption~\ref{ass:localization},
\[
\varphi_t=\Prob(\Eset_t^{\mathrm c})\longrightarrow0
\qquad\text{as }t\downarrow0.
\]
\end{lemma}
\begin{proof}[Proof of Lemma~\ref{lem:energy}]
Using \eqref{eq:global-time-change},
\[
\min_i|X_t-a_ty_i|
=a_t\dist\left(\frac{Z_{q(t)}}{q(t)},S\right).
\]
By \eqref{eq:global-distance}, this distance is at most
$a_t\Lambda/8$ for every sufficiently small $t$, almost surely.
Thus $\one_{\Eset_t^{\mathrm c}}\to0$ almost surely, and dominated
convergence gives $\Prob(\Eset_t^{\mathrm c})\to0$.
\end{proof}

\subsection{Local propagation and stochastic convolution estimates}\label{sec:localization}

We first bound the stochastic convolution, then use it to control the
trajectory before it exits a moving ball. Fix a deterministic $0<b<T$ and define
\begin{equation}\label{eq:Ndef}
N_b(t)\triangleq\sqrt2\sinh t\int_{T-b}^{T-t}\frac{\dd B_s}{\sinh(T-s)},\qquad0<t\le b.
\end{equation}

For $\eta\in(0,1)$ use the envelopes in \eqref{eq:error-response} at
$\alpha=2\eta\in(0,2)$; the second subscript of $\cN$ is only a confidence
budget, so any value in $(0,2)$ is admissible, and $\eta=\alpha/2$ recovers
Theorem~\ref{thm:main}. Define
\begin{equation}\label{eq:noise-event}
\Gset_{b,\eta}\triangleq\{|N_b(t)|\le\cN_{b,2\eta,d}(t)\text{ for every }0<t\le b\}.
\end{equation}
\begin{proposition}[Uniform bound for the stochastic convolution]\label{prop:envelope}
For every deterministic $0<b<T$ and $\eta\in(0,1)$,
$\Gset_{b,\eta}\in\F_{T-}$ is independent of $\F_{T-b}$, and
\begin{equation}\label{eq:noise-prob}
\Prob(\Gset_{b,\eta}^{\mathrm c}\mid\F_{T-b})
=\Prob(\Gset_{b,\eta}^{\mathrm c})\le\eta\qquad\Prob\text{-a.s.}
\end{equation}
\end{proposition}
\begin{proof}[Proof of Proposition~\ref{prop:envelope}]
Set
\begin{equation}\label{eq:Mb}
M_b(r)=\int_{T-b}^r\frac{\dd B_s}{\sinh(T-s)},\qquad T-b\le r<T.
\end{equation}
Its deterministic covariations are
\begin{equation}\label{eq:quadratic-variation}
\langle M_b^{(j)},M_b^{(k)}\rangle_{T-t}
=\delta_{jk}\vartheta_b(t),\qquad
\vartheta_b(t)\triangleq\int_t^b\sinh^{-2}\sigma\dd \sigma=\coth t-\coth b.
\end{equation}
By L\'evy's characterization, the deterministic time change gives a standard
Brownian motion $\Bset$, independent of $\F_{T-b}$, such that almost surely,
simultaneously for $0<t\le b$,
\begin{equation}\label{eq:N-representation}
M_b(T-t)=\Bset_{\vartheta_b(t)},\qquad
N_b(t)=\sqrt2\sinh t\,\Bset_{\vartheta_b(t)}.
\end{equation}
Independence follows because $M_b$ uses only Brownian increments after
$T-b$ with a deterministic integrand.

For $V,L>0$, $\exp(|\Bset_{\lambda}|^2/(4V))$, $0\le \lambda\le V$, is a nonnegative
submartingale with terminal expectation $2^{d/2}$. Doob's inequality gives
\begin{equation}\label{eq:max-bound}
\Prob\left(\sup_{\lambda\le V}|\Bset_{\lambda}|\ge2\sqrt{VL}\right)
\le2^{d/2}e^{-L}.
\end{equation}
Recall that $\ell_{2\eta,d}$ is defined in \eqref{eq:error-response} with
$\alpha=2\eta$. For $k\ge0$, put
\begin{equation}\label{eq:segment-data}
\begin{aligned}
t_k&=b2^{-k},&\eta_k&=\frac{6\eta}{\pi^2(k+1)^2},\\
L_k&=\ell_{2\eta,d}+2\log(k+1),&V_k&=\vartheta_b(t_{k+1}).
\end{aligned}
\end{equation}
Then $2^{d/2}e^{-L_k}=\eta_k$, and
\begin{equation}\label{eq:segment}
\mathcal S_k=\{\sup_{\lambda\le V_k}|\Bset_{\lambda}|\le2\sqrt{V_kL_k}\},
\qquad \Prob(\mathcal S_k^{\mathrm c})\le\eta_k.
\end{equation}
For $t\in[t_{k+1},t_k]$, $\vartheta_b(t)\le V_k\le2/t$,
$\sinh t\le c(b)t$, and
$L_k\le\ell_{2\eta,d}+2\log(1+\log_2(b/t))$. Consequently, on
$\mathcal S_k$,
\[
|N_b(t)|\le2\sqrt2c(b)t\sqrt{V_kL_k}
\le4c(b)\sqrt{tL_k}\le\cN_{b,2\eta,d}(t).
\]
The union bound and $\sum_k\eta_k=\eta$ prove the probability estimate.
\end{proof}

\Needspace{6\baselineskip}
The following proposition is the main result of this subsection.

\begin{proposition}[Local control before exit]\label{prop:local-control}
Under Assumption~\ref{ass:localization}\textup{(\ref{ass:local})}, fix a deterministic $0<b<T$, and put
$\rho(\sigma)=a_\sigma\Lambda/4$. There is a full-probability event $\Omega_b^0\in\F_{T-}$ and,
for each $\eta\in(0,1)$, an event $\Gset_{b,\eta}\in\F_{T-}$ independent of
$\F_{T-b}$ with $\Prob(\Gset_{b,\eta}\mid\F_{T-b})=\Prob(\Gset_{b,\eta})\ge1-\eta$.
On $\Omega_b^0\cap\Gset_{b,\eta}$, simultaneously for every $i$ and $0<t\le b$,
\[
|\widehat X_{T-b}-a_by_i|\le\frac{a_b\Lambda}{8},\qquad
|X_{\sigma}-a_{\sigma}y_i|\le\rho(\sigma)\quad(\sigma\in[t,b])
\]
imply
\begin{equation}\label{eq:local-common-envelope}
|X_t-a_ty_i|\le G_b(t)+R_E(t,b)+\cN_{b,2\eta,d}(t).
\end{equation}
The envelopes are defined in \eqref{eq:Gb}--\eqref{eq:error-response}, using the
same formulas at $\alpha=2\eta\in(0,2)$; $\Gset_{b,\eta}$ is defined in
\eqref{eq:noise-event}.
\end{proposition}

\begin{proof}[Proof of Proposition~\ref{prop:local-control}]
We first derive a propagation identity for the displacement from a moving center.
Under Assumption~\ref{ass:localization}\textup{(\ref{ass:local})}, fix $0<b<T$ and put
$Z_i(r)=\widehat X_r-a_{T-r}y_i$.
Writing $t=T-r$, the identities
$\frac{\dd}{\dd r}(a_{T-r}y_i)=a_ty_i$ and $1-2/h_t=-\coth t$ give
\[
\dd Z_i=\left[-\coth t\,Z_i+
\frac{2(\mu(\widehat X_r,t)-a_ty_i)}{h_t}+2e_r\right]\dd r+\sqrt2\,\dd B_r.
\]
By variation of constants,
\[
\dd\left(\frac{Z_i(r)}{\sinh(T-r)}\right)
=\frac{2(\mu(\widehat X_r,T-r)-a_{T-r}y_i)}{h_{T-r}\sinh(T-r)}\dd r
+\frac{2e_r}{\sinh(T-r)}\dd r
+\frac{\sqrt2}{\sinh(T-r)}\dd B_r.
\]
Integrating from $T-b$ to $T-t$, multiplying by $\sinh t$ and substituting
$\sigma=T-s$ in the Lebesgue integrals, we obtain, on a single full-probability
event $\Omega_b^0$, simultaneously for every $i\in\{1,\ldots,N\}$ and $0<t\le b$,
\begin{equation}\label{eq:propagation}
\begin{aligned}
Z_i(T-t)={}&\frac{\sinh t}{\sinh b}Z_i(T-b)
+2\sinh t\int_t^b\frac{\mu(X_{\sigma},\sigma)-a_{\sigma}y_i}{h_\sigma\sinh \sigma}\dd \sigma\\
&+2\sinh t\int_t^b\frac{e_{T-\sigma}}{\sinh \sigma}\dd \sigma+N_b(t).
\end{aligned}
\end{equation}
Here $N_b$ is defined in \eqref{eq:Ndef}.

We estimate the geometric term in \eqref{eq:propagation} while the path
remains in the moving ball.
For $i\ne j$, put $d_{ij}=|y_i-y_j|$,
$A_{ij}=\tfrac12d_{ij}(d_{ij}-\Lambda/2)$ and $\chi(\sigma)=(e^{2\sigma}-1)^{-1}$.
For $\sigma>0$ and $z=x-a_{\sigma}y_i$ with $|z|\le a_\sigma\Lambda/4$,
\[
|x-a_{\sigma}y_j|^2-|z|^2
=a_\sigma^2d_{ij}^2-2a_{\sigma}z\cdot(y_j-y_i)
\ge a_\sigma^2d_{ij}(d_{ij}-\Lambda/2)=2a_\sigma^2A_{ij}.
\]
Since $a_\sigma^2/h_\sigma=\chi(\sigma)$ and $\omega_i\le1$, the weight ratio and the identity
\[
\mu(x,\sigma)-a_{\sigma}y_i
=\sum_{j\ne i}\omega_j(x,\sigma) a_\sigma(y_j-y_i)
\]
yield
\begin{equation}\label{eq:weights-ratio}
\frac{\omega_j(x,\sigma)}{\omega_i(x,\sigma)}\le e^{-A_{ij}\chi(\sigma)}\ (j\ne i),\qquad
|\mu(x,\sigma)-a_{\sigma}y_i|\le F_i(\sigma)\triangleq a_\sigma\sum_{j\ne i}d_{ij}e^{-A_{ij}\chi(\sigma)}.
\end{equation}
Thus, if the path remains in the moving ball throughout $[t,b]$, the norm
of the second term on the right-hand side of \eqref{eq:propagation} is
bounded by
\begin{equation}\label{eq:geom-response}
R_i^{\geom}(t,b)\triangleq2\sinh t\int_t^b\frac{F_i(\sigma)}{h_\sigma\sinh \sigma}\dd \sigma.
\end{equation}
Using $h_\sigma=2a_\sigma\sinh \sigma$ and $\chi'(\sigma)=-1/(2\sinh^2\sigma)$, direct integration gives
\begin{equation}\label{eq:geom-closed}
R_i^{\geom}(t,b)=2\sinh t\sum_{j\ne i}\frac{d_{ij}}{A_{ij}}
\left[e^{-A_{ij}\chi(b)}-e^{-A_{ij}\chi(t)}\right].
\end{equation}
Since $d_{ij}\ge\Lambda$, we have $A_{ij}\ge\Lambda^2/4$ and
$d_{ij}/A_{ij}\le4/\Lambda$. Hence
\[
2\sum_{j\ne i}\frac{d_{ij}}{A_{ij}}e^{-A_{ij}\chi(b)}
\le\frac{8(N-1)}{\Lambda}e^{-\Lambda^2\chi(b)/4}=g(b).
\]
Together with $\chi(t)\ge\chi(b)$, this yields
\begin{equation}\label{eq:geom-coarse-proof}
0\le R_i^{\geom}(t,b)\le\sinh(t)g(b).
\end{equation}

For the estimate of the score-error term, the envelope \eqref{eq:envelope}
shows that, on $\Omega_b^0$, the third term on the
right-hand side of \eqref{eq:propagation} satisfies, simultaneously for $0<t\le b$,
\[
\left|2\sinh t\int_t^b\frac{e_{T-\sigma}}{\sinh \sigma}\dd \sigma\right|
\le2\sinh t\int_t^b\frac{E(\sigma)}{\sinh \sigma}\dd \sigma=R_E(t,b),
\]
where $R_E$ is defined in \eqref{eq:error-response}.

On $\Omega_b^0\cap\Gset_{b,\eta}$, combine \eqref{eq:propagation},
the preceding bounds, and Proposition~\ref{prop:envelope}.
Whenever the moving-ball condition holds throughout $[t,b]$, this gives,
simultaneously for every $i$ and $0<t\le b$,
\begin{equation}\label{eq:local-control}
\begin{aligned}
|X_t-a_ty_i|&\le\frac{\sinh t}{\sinh b}|\widehat X_{T-b}-a_by_i|
+R_i^{\geom}(t,b)\\
&\quad+R_E(t,b)+\cN_{b,2\eta,d}(t).
\end{aligned}
\end{equation}
The independence and probability claims follow from \eqref{eq:noise-prob}.
Under the entrance condition, the first two terms are at most $G_b(t)$ by \eqref{eq:Gb} and
\eqref{eq:geom-coarse-proof}, proving \eqref{eq:local-common-envelope}.
\end{proof}

\subsection{Completion of the proof}\label{sec:proof-main}

The local estimate assumes that the trajectory stays inside a moving ball.
After choosing the terminal window sufficiently small, a first exit would
contradict this estimate: at the exit time, it would place the trajectory
strictly inside the ball. We make this argument precise below.

\begin{lemma}[Choice of a terminal window]\label{lem:window-selection}
Under Assumption~\ref{ass:localization}, for every
$\eta\in(0,1)$ and $b_0>0$, there exists a deterministic
$0<b<\min\{T,b_0\}$ such that
\begin{equation}\label{eq:conditions}
\begin{gathered}
\sinh(b)g(b)\le\frac{a_b\Lambda}{32},\qquad
\sup_{0<t\le b}R_E(t,b)\le\frac{a_b\Lambda}{32},\\
\sup_{0<t\le b}\cN_{b,2\eta,d}(t)\le\frac{a_b\Lambda}{32},\qquad
\varphi_b\le\eta.
\end{gathered}
\end{equation}
\end{lemma}
\begin{proof}[Proof of Lemma~\ref{lem:window-selection}]
Since \(1\le\sinh \sigma/\sigma\le c(b)\) on \((0,b]\), the weighted response
satisfies
\begin{equation}\label{eq:error-weighted}
\begin{aligned}
 \frac{2}{c(b)}t\int_t^b\frac{E(\sigma)}{\sigma}\dd \sigma&\le R_E(t,b)\le2c(b)t\int_t^b\frac{E(\sigma)}{\sigma}\dd \sigma,\\
 \sup_{0<t\le b}R_E(t,b)&\le2c(b)\sup_{0<t\le b}
 \left[t\int_t^T\frac{E(\sigma)}{\sigma}\dd \sigma\right].
\end{aligned}
\end{equation}
The weighted condition~\eqref{eq:weighted-vanish} therefore gives
\begin{equation}\label{eq:error-decay}
 \sup_{0<t\le b}R_E(t,b)\longrightarrow0\qquad(b\downarrow0).
\end{equation}
For the noise envelope, set $\lambda=\log_2(b/t)\ge0$. Then
\[
\frac{\cN_{b,2\eta,d}(t)^2}{16c(b)^2}
=b2^{-\lambda}[\ell_{2\eta,d}+2\log(1+\lambda)]
\le b2^{-\lambda}(\ell_{2\eta,d}+2\lambda)
\le b\left(\ell_{2\eta,d}+\frac{2}{\mathrm e\log2}\right),
\]
since $2\lambda2^{-\lambda}\le2/(\mathrm e\log2)$. Thus
\[
\sup_{0<t\le b}\cN_{b,2\eta,d}(t)
\le4c(b)\sqrt{b\left(\ell_{2\eta,d}+\frac{2}{\mathrm e\log2}\right)}.
\]
Set $b_k=\min\{T/2,1\}2^{-k}$. The first three left-hand sides in
\eqref{eq:conditions} tend to zero, whereas $a_{b_k}\Lambda/32\to\Lambda/32>0$.
Lemma~\ref{lem:energy} gives \(\varphi_{b_k}\to0\), so the final
condition holds for all sufficiently large \(k\).
Choose such a $k$ with $b_k<b_0$; then all estimates in
\eqref{eq:conditions} hold.
\end{proof}

\begin{proof}[Proof of Theorem~\ref{thm:main}]
\stepheading{Window selection.}\label{lem:common}
Fix $\alpha\in(0,1)$.
Apply Lemma~\ref{lem:window-selection} with $\eta=\alpha/2$ and $b_0=T$ to choose
$b\in(0,T)$ satisfying \eqref{eq:conditions}.
Then $\varphi_b\le\eta$.

\proofstep{No exit.}\label{prop:noexit}
For the chosen window $b$, recall $I_b,\rho_b^{\mathrm{ent}},\Eset_b$ from
\eqref{eq:entry}--\eqref{eq:entrance-event} and put
$\rho(t)=a_t\Lambda/4$ and $Z_i(r)=\widehat X_r-a_{T-r}y_i$. Define
\begin{equation}\label{eq:exit-times}
\sigma_i^b=\inf\{r\in[T-b,T):|Z_i(r)|\ge\rho(T-r)\},\qquad\inf\varnothing=T,
\qquad\sigma_{I_b}^b=\sigma_i^b\ \text{on }\{I_b=i\}.
\end{equation}
The argument in this step uses only the entrance event and the first three
local bounds in \eqref{eq:conditions}; it does not use the global weighted
condition. Recall $\rho_b^{\mathrm{ent}}=a_b\Lambda/8$ from \eqref{eq:entry} and
$\cR_{b,2\eta}$ from \eqref{eq:main-envelope}. By \eqref{eq:Gb} and the first three estimates in
\eqref{eq:conditions}, for all $0<t\le b$,
\[
\begin{aligned}
\cR_{b,2\eta}(t)
&=G_b(t)+R_E(t,b)+\cN_{b,2\eta,d}(t)\\
&=\frac{\sinh t}{\sinh b}\rho_b^{\mathrm{ent}}+\sinh t\,g(b)
+R_E(t,b)+\cN_{b,2\eta,d}(t)
\le\frac{7a_b\Lambda}{32}.
\end{aligned}
\]

Work on $\Omega_b^0\cap\Eset_b\cap\Gset_{b,\eta}$, with the same
full-probability event $\Omega_b^0$ as in
Proposition~\ref{prop:local-control}, and put $i\triangleq I_b$.
Suppose $\sigma\triangleq \sigma_i^b<T$ and set $t_\sigma\triangleq T-\sigma$.
Since $|Z_i(T-b)|\le \rho_b^{\mathrm{ent}}<\rho(b)$ on $\Eset_b$, continuity gives
\[
|Z_i(\sigma)|=\rho(t_\sigma)
=\frac{a_{t_\sigma}\Lambda}{4}
\ge\frac{a_b\Lambda}{4}=\frac{8a_b\Lambda}{32}.
\]
On $[T-b,\sigma]$, the first-exit definition gives
$|Z_i(r)|\le\rho(T-r)$. Proposition~\ref{prop:local-control} applies
simultaneously to all indices and times, so its estimate
\eqref{eq:local-common-envelope} at $t=t_\sigma$ yields
\begin{equation}\label{eq:exit-contradiction}
\begin{aligned}
|Z_i(\sigma)|
&\le G_b(t_\sigma)+R_E(t_\sigma,b)+\cN_{b,2\eta,d}(t_\sigma)\\
&=\cR_{b,2\eta}(t_\sigma)\le\frac{7a_b\Lambda}{32}.
\end{aligned}
\end{equation}
This contradicts the boundary value $|Z_i(\sigma)|\ge8a_b\Lambda/32$.
Therefore $\sigma_{I_b}^b=T$. Applying Proposition~\ref{prop:local-control}
on every $[t,b]$ now gives
\begin{equation}\label{eq:noexit-envelope}
    |X_t-a_ty_{I_b}|
    \le \cR_{b,2\eta}(t)
    \le\frac{7a_b\Lambda}{32}<\rho(t),
    \qquad 0<t\le b.
\end{equation}
For $j\ne I_b$, the triangle inequality gives
\[
|X_t-a_ty_j|
\ge a_t\Lambda-|X_t-a_ty_{I_b}|
>\frac{3a_t\Lambda}{4}>|X_t-a_ty_{I_b}|,
\]
so $I_b$ is the index of the unique nearest moving center throughout the window,
as asserted in part \textup{(i)} of the theorem.

\proofstep{Localization event and terminal limit.}
With $\eta=\alpha/2$, set
\begin{equation}\label{eq:main-event}
\Hset_{b,\alpha}\triangleq\Omega_b^0\cap\Eset_b\cap\Gset_{b,\eta}\in\F_{T-}.
\end{equation}
Since $\Prob(\Omega_b^0)=1$, $\Eset_b\in\F_{T-b}$, and
$\Gset_{b,\eta}$ is independent of $\F_{T-b}$, window selection gives
\begin{equation}\label{eq:main-prob}
\begin{aligned}
\Prob(\Hset_{b,\alpha})
&=\Prob(\Eset_b)\Prob(\Gset_{b,\eta})\\
&\ge(1-\varphi_b)(1-\eta)\ge(1-\eta)^2\ge1-\alpha.
\end{aligned}
\end{equation}
Step~\ref{prop:noexit} gives no exit and the uniform estimate on this event.

For fixed $b$, \eqref{eq:error-weighted} and
the weighted condition~\eqref{eq:weighted-vanish} give
\begin{equation}\label{eq:error-terminal-decay}
 R_E(t,b)\le2c(b)t\int_t^b\frac{E(\sigma)}{\sigma}\dd \sigma\longrightarrow0\qquad(t\downarrow0).
\end{equation}
Also, $G_b(t)\to0$ directly, and
$\cN_{b,\alpha,d}(t)\to0$ by its definition in \eqref{eq:error-response}.
All three terms are continuous on $(0,b]$, so $\cR_{b,\alpha}(t)\to0$ and
\[
|X_t-y_{I_b}|
\le \cR_{b,\alpha}(t)+(1-a_t)R\to0.
\]

\proofstep{Identification of the terminal index.}
Proposition~\ref{prop:global-localization} already gives an
\(\F_{T-}\)-measurable index \(I_\infty\) with
\(X_t\to y_{I_\infty}\) almost surely as \(t\downarrow0\). The limit established in Step~3
and distinctness of the training points imply
\(I_b=I_\infty\) almost surely on \(\Hset_{b,\alpha}\).
On the probability-one convergence event, both
\(\lvert X_t-y_{I_\infty}\rvert\) and
\(\lvert X_t-a_ty_{I_\infty}\rvert\) tend to zero. Their corresponding
center separations are \(\Lambda\) and \(a_t\Lambda\to\Lambda\), so the
triangle inequality makes every other center strictly farther away for all
sufficiently small \(t\). This proves part~\textup{(ii)}.
The definition \(X_0=y_{I_\infty}\) gives an almost-sure continuous
extension to \([0,T]\).
\end{proof}

\subsection{Positivity of terminal selection probabilities}
\label{sec:positive-terminal-support}

We combine the exact diffusion's positive transition density with a direct
application of Girsanov's theorem on a truncated interval and the local
no-exit estimate.

\begin{proof}[Proof of Corollary~\ref{cor:positive-terminal-support}]
Apply Lemma~\ref{lem:window-selection} with $\eta=1/4$ to choose a deterministic
$0<b<T$ satisfying the first three conditions in \eqref{eq:conditions},
and put $H=T-b$.
For each $i$, put
\[
D_i\triangleq B(a_by_i,\rho_b^{\mathrm{ent}}),\qquad
\rho_b^{\mathrm{ent}}=\frac{a_b\Lambda}{8},\qquad
\mathcal A_{b,i}^{\mathrm{ent}}\triangleq\{\widehat X_{T-b}\in D_i\}.
\]
The center separation is $a_b\Lambda=8\rho_b^{\mathrm{ent}}$, so
\[
\mathcal A_{b,i}^{\mathrm{ent}}\in\F_{T-b},\qquad
\mathcal A_{b,i}^{\mathrm{ent}}\subset\Eset_b\cap\{I_b=i\}.
\]

Let $P_b^x$ be the path law on $C([0,H];\R^d)$ of the exact equation
\begin{equation}\label{eq:exact-sde}
    \dd Z_r=[Z_r+2s(Z_r,T-r)]\dd r+\sqrt2\,\dd B_r,
    \qquad 0\le r\le H,\qquad Z_0=x,
\end{equation}
and write $P_b^{v_T}\triangleq\int P_b^x\,v_T(\dd x)$.
Well-posedness and weak uniqueness follow from Appendix~\ref{app:wellposed}.
Time reversal of the forward Ornstein--Uhlenbeck process and weak uniqueness
identify the exact diffusion started from $w(T)$ with
$(X^{\mathrm f}_{T-r})_{0\le r\le H}$ in law.
The joint density of $(X_b^{\mathrm f},X_T^{\mathrm f})$ is
$w(z,b)\varphi_{h_{T-b}}(x-a_{T-b}z)$. Thus time reversal and Bayes' formula give,
for $w(T)$-almost every $x$ and every Borel set $D\subset\R^d$,
\begin{equation}\label{eq:reverse-transition-identity}
P_b^x\{\omega(T-b)\in D\}
=\frac{\int_D w(z,b)\varphi_{h_{T-b}}(x-a_{T-b}z)\dd z}{w(x,T)}
=\int_D k_{0,T-b}(x,z)\dd z,
\end{equation}
where
\begin{equation}\label{eq:reverse-transition-kernel}
k_{0,T-b}(x,z)\triangleq
\frac{w(z,b)\varphi_{h_{T-b}}(x-a_{T-b}z)}{w(x,T)}>0,
\qquad \int_{\R^d}k_{0,T-b}(x,z)\dd z=1.
\end{equation}
In fact, \eqref{eq:reverse-transition-identity} holds for every $x$, by
continuous dependence on the initial condition. To see this, let $Z^x$ solve
\eqref{eq:exact-sde} from $x$ and couple $Z^x,Z^{x'}$ with the same Brownian
motion. By \eqref{eq:lipschitz} and Gronwall's inequality,
\[
\sup_{0\le r\le T-b}|Z_r^x-Z_r^{x'}|
\le e^{L_b(T-b)}|x-x'|.
\]
For every continuous compactly supported function $f$, the expectation
$\E[f(Z_{T-b}^x)]$ is continuous in $x$.
Fix any $x\in\R^d$. Since $w(\cdot,T)>0$, the set of initial points where
\eqref{eq:reverse-transition-identity} holds is dense in $\R^d$. We can therefore
choose $x_n\to x$ from this set. Passing to the limit, using the explicit
formula for $k_{0,T-b}$, gives
\[
\begin{aligned}
\E[f(Z_{T-b}^x)]
&=\lim_{n\to\infty}\E[f(Z_{T-b}^{x_n})]\\
&=\lim_{n\to\infty}\int_{\R^d}f(z)k_{0,T-b}(x_n,z)\dd z\\
&=\int_{\R^d}f(z)k_{0,T-b}(x,z)\dd z.
\end{aligned}
\]
Thus the law of $Z_{T-b}^x$ and the measure $k_{0,T-b}(x,z)\dd z$ give the
same integral for every continuous compactly supported $f$, so they are
equal. Consequently, \eqref{eq:reverse-transition-identity} holds for every
$x\in\R^d$ and every Borel set $D\subset\R^d$.
Since $D_i$ has positive Lebesgue measure, it follows that
\[
P_b^{v_T}\{\omega(T-b)\in D_i\}
=\int_{\R^d}v_T(\dd x)\int_{D_i}k_{0,T-b}(x,z)\dd z>0.
\]
To transfer this positive probability to the perturbed process, note that
the deterministic envelope gives Novikov's condition:
\[
\E_{\Prob}\exp\left(\int_0^H|e_r|^2\dd r\right)
\le\exp\left(\int_b^T E(\sigma)^2\dd\sigma\right)<\infty.
\]
Consequently,
\[
M_r\triangleq\exp\left\{-\sqrt2\int_0^r e_\ell\cdot\dd B_\ell
-\int_0^r|e_\ell|^2\dd\ell\right\},\qquad 0\le r\le H,
\]
is a positive martingale with $M_0=1$.
Define $\Prob^M$ on $\F_H$ by $\dd\Prob^M=M_H\dd\Prob$.
Girsanov's theorem~\cite{girsanov1960} makes
$B_r^M=B_r+\sqrt2\int_0^r e_\ell\dd\ell$ Brownian under $\Prob^M$,
and $\widehat X$ then solves the exact equation~\eqref{eq:exact-sde}.
Moreover, $\E_{\Prob}[M_H\mid\F_0]=1$ preserves the initialization law $v_T$.
Weak uniqueness therefore identifies the path law of $\widehat X$ under
$\Prob^M$ with $P_b^{v_T}$, so
\[
\Prob^M(\mathcal A_{b,i}^{\mathrm{ent}})
=P_b^{v_T}\{\omega(H)\in D_i\}>0.
\]
Since $\Prob^M\ll\Prob|_{\F_H}$, this implies
$\Prob(\mathcal A_{b,i}^{\mathrm{ent}})>0$.

The no-exit argument in the proof of Theorem~\ref{thm:main}, using the first
three bounds in \eqref{eq:conditions}, shows that on
$\Omega_b^0\cap\mathcal A_{b,i}^{\mathrm{ent}}\cap\Gset_{b,\eta}$ the trajectory
stays in the moving ball about $a_ty_i$ for all $0<t\le b$ and converges to
$y_i$. Here $\Prob(\Omega_b^0)=1$, and $\Gset_{b,\eta}$ is independent of
$\F_{T-b}$ with probability at least $1-\eta=3/4$. Since
$\mathcal A_{b,i}^{\mathrm{ent}}\in\F_{T-b}$, we obtain
\begin{equation}\label{eq:positive-atom-bound}
\begin{aligned}
q_i=\Prob(I_\infty=i)
&\ge\Prob(\mathcal A_{b,i}^{\mathrm{ent}}\cap\Gset_{b,\eta})\\
&=\Prob(\mathcal A_{b,i}^{\mathrm{ent}})\Prob(\Gset_{b,\eta})
\ge\frac34\Prob(\mathcal A_{b,i}^{\mathrm{ent}})>0.
\end{aligned}
\end{equation}
Finally, $X_0=y_{I_\infty}$ almost surely and distinctness of the training
points give
\[
\Prob(X_0\in\cdot)=\sum_{i=1}^Nq_i\delta_{y_i},\qquad
\operatorname{supp}(\Prob\circ X_0^{-1})=\{y_1,\ldots,y_N\}=S,
\]
which proves \eqref{eq:terminal-law} and the support assertion.
\end{proof}

\subsection{A two-atom example with altered weights}\label{sec:two-atom-proof}

\begin{proof}[Proof of Example~\ref{ex:changed-weights}]
Put $c_p=\frac12\log(p/(1-p))$. The identity
\[
 |\tanh(u+c_p)-\tanh u|
 =\frac{2|\sinh c_p|}{\cosh(2u+c_p)+\cosh c_p}
 \le2\tanh(|c_p|/2)
\]
is an equality at $u=-c_p/2$. Together with $h_t=2a_t\sinh t$,
this proves~\eqref{eq:two-atom-envelope}. The envelope is bounded and
vanishes on $(0,b)$, so~\eqref{eq:weighted-vanish} holds. The smooth
mixture scores, joined at the deterministic switching time, give
compatible strong solutions on every truncated interval.

Let $P$ be the exact reverse law initialized from $w(T)$, with its
natural filtration, and let $Y=\lim_{t\downarrow0}X_t\in\{-1,1\}$.
Bayes' formula for the forward OU process gives
$m_t(X_t)=\E_P[Y\mid\F_{T-t}]$, a bounded martingale in generation time.
Define the full-path law $Q_{p,b}$ by
\begin{equation}\label{eq:two-atom-density}
 \frac{\dd Q_{p,b}}{\dd P}
 =\frac{1+\theta m_b(X_b)}{1+\theta m_T(X_T)}.
\end{equation}
The density and its reciprocal are bounded by positive constants
because $|\theta|<1$. Its conditional expectation given $X_T$ is one,
so the entire initialization law is preserved. Until generation time
$T-b$, its density process is
$L_{T-t}=(1+\theta m_t(X_t))/(1+\theta m_T(X_T))$; afterwards it is
constant. It\^o's formula gives
\[
 \frac{\dd L_r}{L_r}=\sqrt2\,e^{p,b}(\widehat X_r,T-r)\dd B_r.
\]
Girsanov's theorem therefore adds precisely $2e^{p,b}$ to the drift,
and adds nothing on $0<t<b$. This identifies $Q_{p,b}$ with the
stated process. Equivalence with $P$ on full path space gives the
terminal limit and positive probability of each endpoint.

Set $V_b(z)=\Var_P(m_b(X_b)\mid X_T=z)$, using the explicit conditional
kernel
\begin{equation}\label{eq:two-atom-variance}
 \begin{aligned}
 K_b(x\mid z)&=\frac{w(x,b)\varphi_{h_{T-b}}(z-a_{T-b}x)}{w(z,T)},\\
 V_b(z)&=\int_\R m_b(x)^2K_b(x\mid z)\dd x-m_T(z)^2.
 \end{aligned}
\end{equation}
Since $K_b$ is strictly positive and $m_b$ is strictly increasing,
$0<V_b(z)<1-m_T(z)^2$ for every finite $z$. Reflection makes $V_b$ even.
For $U=m_b(X_b)$ and $m=m_T(X_T)$, conditional expectation in
\eqref{eq:two-atom-density} yields
\[
 q_+(p,b)=\E_P\!\left[\frac{(1+U)(1+\theta U)}{2(1+\theta m)}\right]
 =\frac12+\frac\theta2\E_P\!\left[\frac{V_b(X_T)}{1+\theta m_T(X_T)}\right].
\]
Averaging the last integrand at $z$ and $-z$ gives the exact formula
\begin{equation}\label{eq:two-atom-weight}
 q_+(p,b)=\frac12+\frac\theta2\int_\R w(z,T)
 \frac{V_b(z)}{1-\theta^2m_T(z)^2}\dd z.
\end{equation}
This gives $q_+(1/2,b)=1/2$ and $q_+(p,b)>1/2$ for $p>1/2$.
Differentiation under the integral, justified locally uniformly for
$|\theta|<1$, gives
\[
 \partial_pq_+(p,b)=\int_\R w(z,T)V_b(z)
 \frac{1+\theta^2m_T(z)^2}{(1-\theta^2m_T(z)^2)^2}\dd z>0.
\]
For numerical evaluation of~\eqref{eq:two-atom-variance}, let
$\rho=a_{T-b}$, $h_g=h_{T-b}$,
$\beta=\rho h_b/h_T$, $c=a_bh_g/h_T$ and $\nu=h_bh_g/h_T>0$.
The same kernel is the Gaussian bridge mixture
\[
 K_b(x\mid z)=\frac{1+m_T(z)}2\varphi_\nu(x-\beta z-c)
             +\frac{1-m_T(z)}2\varphi_\nu(x-\beta z+c).
\]
Thus~\eqref{eq:two-atom-weight} is an explicit Gaussian quadrature
formula for the changed selection probability.
\end{proof}

\section{Power-law asymptotics and critical displacement}\label{sec:rates}

For the power-law perturbation \eqref{eq:power-error} with $1/2<\gamma<1$,
we compute the leading term of the displacement from the terminal atom.
For the critical perturbation \eqref{eq:critical-error}, the same calculation
yields a nonzero limiting displacement.

For a fixed-direction error
\[
 e_{T-\sigma}=M\sigma^{-\gamma}\xi,\qquad M>0,\quad |\xi|=1,\quad \gamma>0,
\]
using the small-time asymptotic form of the contraction kernel \(t/\sigma \) suggests the response
\[
 2t\int_t^b\frac{M\sigma^{-\gamma}}{\sigma}\dd \sigma \,\xi
 =\frac{2M}{\gamma}\bigl(t^{1-\gamma}-tb^{-\gamma}\bigr)\xi,
 \qquad 0<t\le b.
\]
The proofs below retain the exact kernel $\sinh t/\sinh\sigma$. This asymptotic response vanishes when \(0<\gamma<1\), whereas at \(\gamma=1\) a
persistent fixed-direction perturbation can leave a nonzero limiting
response. Proposition~\ref{prop:matching} identifies the corresponding
leading displacement in the stochastic system when \(1/2<\gamma<1\).

\subsection{Proof of Proposition~\ref{prop:matching}: power-law asymptotics}\label{sec:error-response}

\begin{proof}[Proof of Proposition~\ref{prop:matching}]
For $E(\sigma)=M\sigma^{-\gamma}$ and $1/2<\gamma<1$,
\[
 t\int_t^T\frac{E(\sigma)}{\sigma}\dd \sigma =\frac{M}{\gamma}
 \bigl(t^{1-\gamma}-tT^{-\gamma}\bigr)\longrightarrow0.
\]
Thus Theorem~\ref{thm:main} gives $X_t\to y_{I_\infty}$ almost surely as $t\downarrow0$.

\proofstep{asymptotics of the error term.}
The perturbation acts in the fixed direction $\xi$, so it suffices to compute
its scalar coefficient in \eqref{eq:propagation}.
For fixed $0<b<T$, put $\psi(\sigma)=\sigma /\sinh \sigma $ and $\psi(0)=1$.
Since $\psi(\sigma)=1+O(\sigma^2)$,
\[
\begin{aligned}
H_{M,\gamma}(t,b)&\triangleq
2M\sinh t\int_t^b\frac{\sigma^{-\gamma}}{\sinh \sigma}\dd \sigma \\
&=2M\sinh t\left[
\frac{t^{-\gamma}-b^{-\gamma}}{\gamma}
+\int_t^b \sigma^{-\gamma-1}(\psi(\sigma)-1)\dd \sigma \right].
\end{aligned}
\]
The last integral is bounded by $C_bb^{2-\gamma}/(2-\gamma)$.
Using $\sinh t=t+O(t^3)$ gives
\begin{equation}\label{eq:H-expansion}
H_{M,\gamma}(t,b)=\frac{2M}{\gamma}t^{1-\gamma}+O_b(t).
\end{equation}
Here the subscript $b$ in $O_b(t)$ indicates that the implicit constant
may depend on the fixed $b$, but not on $t$.

\proofstep{choice of terminal windows.}
We apply this expansion on localization events for the fixed process driven
by \eqref{eq:power-error}.
For $\alpha_m=2^{-m}$ and $\eta_m=\alpha_m/2$, choose $b_m$ by
Step~\ref{lem:common} of the proof of Theorem~\ref{thm:main}.
Intersect the propagation events for these
countably many windows with the terminal convergence event to obtain
a full-probability event $\Omega_*$, and set
\begin{equation}\label{eq:Km}
 K_m\triangleq\Omega_*\cap\Eset_{b_m}\cap\Gset_{b_m,\eta_m}.
\end{equation}
The entrance estimate and noise independence give
$\Prob(K_m)\ge(1-\eta_m)(1-\varphi_{b_m})\ge1-\alpha_m$.
Hence $\Prob(\bigcup_mK_m)=1$. The no-exit argument gives
$K_m\subset\Hset_{b_m,\alpha_m}$ and $I_{b_m}=I_\infty$ on $K_m$.

\proofstep{comparison of the terms.}
It remains to show that the other three terms in the propagation formula
vanish after rescaling.
Fix a path in $K_m$ and write $b=b_m$, $\eta=\eta_m$ and
\[
 V_i(t,b)\triangleq2\sinh t\int_t^b
 \frac{\mu(X_{\sigma},\sigma)-a_{\sigma}y_i}{h_{\sigma}\sinh \sigma}\dd \sigma .
\]
The propagation identity \eqref{eq:propagation} gives
\[
X_t-a_ty_{I_\infty}
=\frac{\sinh t}{\sinh b}Z_{I_\infty}(T-b)
+V_{I_\infty}(t,b)+H_{M,\gamma}(t,b)\xi+N_b(t).
\]
With
\begin{equation}\label{eq:Psi}
\Psi_{b,\eta,d}(t)\triangleq\ell_{2\eta,d}
+2\log\bigl(1+\log_2(b/t)\bigr),
\end{equation}
the entrance, geometric and noise bounds yield
\begin{equation}\label{eq:noise-ratio}
\begin{aligned}
t^{\gamma-1}\left|\frac{\sinh t}{\sinh b}Z_{I_\infty}(T-b)\right|
&\le\frac{c(b)\rho_b^{\mathrm{ent}}}{\sinh b}t^\gamma\longrightarrow0,\\
t^{\gamma-1}|V_{I_\infty}(t,b)|
&\le t^{\gamma-1}R_{I_\infty}^{\geom}(t,b)
\le c(b)g(b)t^\gamma\longrightarrow0,\\
t^{\gamma-1}H_{M,\gamma}(t,b)
&=\frac{2M}{\gamma}+O_b(t^\gamma)\longrightarrow\frac{2M}{\gamma},\\
t^{\gamma-1}|N_b(t)|
&\le4c(b)t^{\gamma-1/2}\sqrt{\Psi_{b,\eta,d}(t)}\longrightarrow0.
\end{aligned}
\end{equation}
The last limit uses $\gamma>1/2$ and
$\Psi_{b,\eta,d}(t)=O_b(1+\log\log(b/t))$.
Consequently,
\begin{equation}\label{eq:Km-limit}
 t^{\gamma-1}(X_t-a_ty_{I_\infty})
 \longrightarrow\frac{2M}{\gamma}\xi.
\end{equation}
The difference between moving and fixed centers contributes at most
$Rt^\gamma\to0$ after scaling. Thus \eqref{eq:matching} holds on
$\bigcup_mK_m$, a probability-one set, which proves the stated vector limit.
\end{proof}

\subsection{Proof of Proposition~\ref{prop:critical}: critical displacement}\label{sec:critical-sec}

For \eqref{eq:critical-error}, set
$E^{(b)}(\sigma)\triangleq\kappa \sigma^{-1}\one_{(0,b]}(\sigma)$. Since
\begin{equation}\label{eq:critical-weighted}
 t\int_t^T\frac{E^{(b)}(\sigma)}{\sigma}\dd \sigma =\kappa(1-t/b)\longrightarrow\kappa>0,
 \qquad0<t\le b.
\end{equation}
condition~\eqref{eq:weighted-vanish} fails for the switched process.
Nevertheless, Proposition~\ref{prop:global-localization} applies to the fixed
unperturbed process with the same initial random variable and Brownian motion.
The processes agree before the switch, so convergence of this unperturbed
process gives the required entrance probability for any initialization law.
\begin{proof}[Proof of Proposition~\ref{prop:critical}]
Fix $\kappa,\alpha,\xi$ as in the statement and put $\eta=\alpha/2$.

\proofstep{high-probability no exit.}
Let $\widehat X^{(0)}$ be the solution with $e\equiv0$, driven by the same
$\widehat X_0$ and $B$, and put $X_t^{(0)}=\widehat X_{T-t}^{(0)}$.
These solutions exist and are pathwise unique on every truncated interval
by Appendix~\ref{app:wellposed}, without assumptions on $v_T$.
For each deterministic $0<b<T$, the error $e^{\kappa,b}$ vanishes on
$[0,T-b)$; its value at the single switch time does not affect the drift
integral. Pathwise uniqueness therefore gives
\begin{equation}\label{eq:critical-coupling}
\widehat X_r^{\kappa,b}=\widehat X_r^{(0)}
\quad\text{for every }0\le r\le T-b\quad\text{a.s.},
\qquad X_b^{\kappa,b}=X_b^{(0)}\quad\text{a.s.}
\end{equation}
Applying Proposition~\ref{prop:global-localization} with $E\equiv0$ to this
single unperturbed process gives $X_t^{(0)}\to y_{I_\infty^{(0)}}$ almost surely
for an index $I_\infty^{(0)}$. Hence
$\min_i|X_b^{(0)}-a_by_i|\to0$ almost surely, while
$a_b\Lambda/8\to\Lambda/8>0$. Dominated convergence yields
\begin{equation}\label{eq:crit-entrance}
\varphi_b^{(0)}\triangleq
\Prob\left(\min_i|X_b^{(0)}-a_by_i|>\frac{a_b\Lambda}{8}\right)
\longrightarrow0\qquad(b\downarrow0).
\end{equation}
By~\eqref{eq:critical-coupling}, the entrance failure probability of
$\widehat X^{\kappa,b}$ at its switch time is exactly $\varphi_b^{(0)}$.
The convergence in~\eqref{eq:crit-entrance} thus concerns one fixed process,
even though the switched processes depend on $b$.

For each candidate window $0<b<T$, \eqref{eq:error-weighted} and the noise
estimate in the proof of Lemma~\ref{lem:window-selection} give
\[
\sup_{0<t\le b}R_{E^{(b)}}(t,b)\le2\kappa c(b),\qquad
\sup_{0<t\le b}\cN_{b,2\eta,d}(t)
\le4c(b)\sqrt{b\left(\ell_{2\eta,d}+\frac{2}{\mathrm e\log2}\right)}.
\]
The geometric and noise bounds tend to zero, while
$2\kappa c(b)\to2\kappa<\Lambda/32$ and $a_b\Lambda/32\to\Lambda/32$.
We can therefore choose a deterministic $b=b(v_T,\kappa,\alpha)$ from
$\min\{T/2,1\}2^{-k}$ such that
\begin{equation}\label{eq:crit-conditions}
\begin{gathered}
\sinh(b)g(b)\le\frac{a_b\Lambda}{32},\qquad
\sup_{0<t\le b}R_{E^{(b)}}(t,b)\le\frac{a_b\Lambda}{32},\\
\sup_{0<t\le b}\cN_{b,2\eta,d}(t)\le\frac{a_b\Lambda}{32},\qquad
\varphi_b^{(0)}\le\eta.
\end{gathered}
\end{equation}
Fix this $b$ and the process $\widehat X^{\kappa,b}$. Put
$\rho(t)=a_t\Lambda/4$, and define $I_b$, $\Eset_b$ and $\sigma_i^b$ for this
process by \eqref{eq:entry}, \eqref{eq:entrance-event} and \eqref{eq:exit-times},
respectively. Let $\Omega_b^0$ be the event in Proposition~\ref{prop:local-control}
and $\Gset_{b,\eta}$ the event in \eqref{eq:noise-event}, and set
\begin{equation}\label{eq:critical-event}
\Hset^{\mathrm{crit}}_{b,\alpha}
\triangleq\Omega_b^0\cap\Eset_b\cap\Gset_{b,\eta}.
\end{equation}
The first-exit argument in Step~\ref{prop:noexit} of the proof of
Theorem~\ref{thm:main}, using Proposition~\ref{prop:local-control} and the
first three bounds in \eqref{eq:crit-conditions}, gives on this event
\[
\sigma_{I_b}^b=T,\qquad
|X^{\kappa,b}_t-a_ty_{I_b}|
\le\frac{7a_b\Lambda}{32}<\rho(t),\qquad0<t\le b.
\]
By the independence of $\Gset_{b,\eta}$ from $\F_{T-b}$,
\[
\Prob(\Hset^{\mathrm{crit}}_{b,\alpha})
=\Prob(\Eset_b)\Prob(\Gset_{b,\eta})
\ge(1-\varphi_b^{(0)})(1-\eta)\ge(1-\eta)^2\ge1-\alpha.
\]

\proofstep{propagation identity on the no-exit event.}
For the fixed process $\widehat X^{\kappa,b}$, put
$Z_i(r)\triangleq \widehat X^{\kappa,b}_r-a_{T-r}y_i$ and define
\begin{equation}\label{eq:critical-mixture-response}
V_i(t,b)\triangleq2\sinh t\int_t^b
\frac{\mu(X^{\kappa,b}_{\sigma},\sigma)-a_{\sigma}y_i}{h_{\sigma}\sinh \sigma}\dd \sigma ,
\qquad0<t\le b,
\end{equation}
and
\begin{equation}\label{eq:H-kappa}
H_\kappa(t,b)\triangleq2\kappa\sinh t\int_t^b\frac{\dd \sigma}{\sigma \sinh \sigma}.
\end{equation}
On $\Hset^{\mathrm{crit}}_{b,\alpha}$, the propagation identity
\eqref{eq:propagation} becomes
\begin{equation}\label{eq:critical-propagation}
X^{\kappa,b}_t-a_ty_{I_b}
=\frac{\sinh t}{\sinh b}Z_{I_b}(T-b)
+V_{I_b}(t,b)+H_\kappa(t,b)\xi+N_b(t).
\end{equation}
Since $|Z_{I_b}(T-b)|\le \rho_b^{\mathrm{ent}}$ on $\Eset_b$,
\[
\left|\frac{\sinh t}{\sinh b}Z_{I_b}(T-b)\right|
\le\frac{\rho_b^{\mathrm{ent}}}{\sinh b}\sinh t\longrightarrow0\qquad(t\downarrow0).
\]
Since the path does not exit,
\eqref{eq:geom-coarse-proof} gives
\[
|V_{I_b}(t,b)|\le R_{I_b}^{\geom}(t,b)\le c(b)g(b)t\longrightarrow0.
\]
For the error response, put $\psi(\sigma)\triangleq \sigma /\sinh \sigma =1+O(\sigma^2)$. Then
\[
\begin{aligned}
H_\kappa(t,b)
&=2\kappa\sinh t\int_t^b\frac{\psi(\sigma)}{\sigma^2}\dd \sigma \\
&=2\kappa\sinh t\left[\frac1t-\frac1b
+\int_t^b\frac{\psi(\sigma)-1}{\sigma^2}\dd \sigma \right].
\end{aligned}
\]
The remainder integral is bounded, so
\begin{equation}\label{eq:H-kappa-limit}
H_\kappa(t,b)=2\kappa+O_b(t)\longrightarrow2\kappa\qquad(t\downarrow0).
\end{equation}

Finally, \eqref{eq:critical-event} gives
$\Hset^{\mathrm{crit}}_{b,\alpha}\subset\Gset_{b,\eta}$. On this event,
\eqref{eq:noise-event} and \eqref{eq:error-response} imply
\[
|N_b(t)|\le\cN_{b,2\eta,d}(t)\longrightarrow0\qquad(t\downarrow0).
\]

\proofstep{terminal displacement.}
Taking $t\downarrow0$ in \eqref{eq:critical-propagation}, using the limits in
Step 2 and $a_t\to1$, we obtain, for almost every path of
$\Hset^{\mathrm{crit}}_{b,\alpha}$,
\[
X^{\kappa,b}_t\longrightarrow y_{I_b}+2\kappa\xi.
\]
Finally,
\[
|y_{I_b}+2\kappa\xi-y_{I_b}|=2\kappa,\qquad
|y_{I_b}+2\kappa\xi-y_j|\ge\Lambda-2\kappa>2\kappa\quad(j\ne I_b),
\]
which proves \eqref{eq:critical-limit}.
\end{proof}


\subsection{Reverse-diffusion simulations for a two-atom target}
\label{sec:numerical}

We simulate \eqref{eq:main-sde} for $d=1$ and $S=\{-1,+1\}$, so that
$u_0=\frac12(\delta_{-1}+\delta_{+1})$ and $\Lambda=2$; by \eqref{eq:score} the
exact score is $\bigl(a_t\tanh(a_tx/h_t)-x\bigr)/h_t$. Every run starts from the
exact initial law $v_T=w(T)$ and uses $T=1$. We compare the persistent errors
$e_r=M(T-r)^{-\gamma}\xi$, with $M=0.02$, $\xi=1$ and
$\gamma\in\{0.25,0.5,0.75\}$, with the switched perturbation
\eqref{eq:critical-error} of Proposition~\ref{prop:critical}, for which
$\kappa=0.01$ and $b=0.05$. The admissibility condition $\kappa<\Lambda/64$
holds, and $2\kappa=0.02$ is the terminal displacement predicted by
\eqref{eq:critical-limit}.

The drift is stiff near $t=0$, so we integrate the transformed system
\eqref{eq:global-time-change}--\eqref{eq:global-transformed-sde}. With
$q=\coth t-1$, $t=t(q)$ and $Z_q=X_t/\sinh t$, it reads
\[
 \dd Z_q=\Bigl[\tanh\frac{Z_q}{2}+2\sinh t\,e_{T-t}\Bigr]\dd q+\sqrt2\,\dd W_q ,
\]
where $\pi_+-\pi_-=\tanh(Z_q/2)$ for this target. On a geometric grid
$q_{k+1}=(1+\varepsilon)q_k$ the bounded drift and the exact Brownian scaling are
resolved uniformly in $1/t$, so that Euler--Maruyama remains stable down to the
stopping time. All runs use $N=10^4$ trajectories, seed $20261028$,
$\varepsilon=5\times10^{-4}$ (29962 steps from $q(1)=0.3130$ to
$q(10^{-6})\approx10^6$) and $t_{\mathrm{stop}}=10^{-6}$; we record
$X_{t_{\mathrm{stop}}}$ and the terminal offset $X_{t_{\mathrm{stop}}}-y_I$,
where $y_I=\operatorname{sign}X_{t_{\mathrm{stop}}}$.

\pgfplotstableset{col sep=comma}

\begin{figure}[htbp]
\centering
\begin{minipage}[t]{0.305\textwidth}
\centering
\begin{tikzpicture}[baseline=(current bounding box.north)]
\begin{semilogxaxis}[
  width=\linewidth, height=0.78\linewidth,
  title={(a) $\gamma=0.25$},
  xlabel={$t$}, ylabel={$X_t$},
  xmin=1e-6, xmax=1, ymin=-2.1, ymax=2.1,
  ytick={-2,-1,0,1,2},
  tick label style={font=\tiny},
  label style={font=\scriptsize},
  title style={font=\small},
  grid=major, major grid style={gray!20},
]
\addplot[dashed, gray, domain=1e-6:1, samples=2, forget plot] {1};
\addplot[dashed, gray, domain=1e-6:1, samples=2, forget plot] {-1};
\addplot[line width=0.35pt, blue!75!black] table[x=t,y=x1] {Python/data/paths_g025.csv};
\addplot[line width=0.35pt, blue!55!black] table[x=t,y=x2] {Python/data/paths_g025.csv};
\addplot[line width=0.35pt, blue!35!black] table[x=t,y=x3] {Python/data/paths_g025.csv};
\addplot[line width=0.35pt, blue!75!black, densely dotted] table[x=t,y=x4] {Python/data/paths_g025.csv};
\addplot[line width=0.35pt, blue!55!black, densely dotted] table[x=t,y=x5] {Python/data/paths_g025.csv};
\addplot[line width=0.35pt, blue!35!black, densely dotted] table[x=t,y=x6] {Python/data/paths_g025.csv};
\end{semilogxaxis}
\end{tikzpicture}
\end{minipage}\hfill
\begin{minipage}[t]{0.305\textwidth}
\centering
\begin{tikzpicture}[baseline=(current bounding box.north)]
\begin{semilogxaxis}[
  width=\linewidth, height=0.78\linewidth,
  title={(b) $\gamma=0.75$},
  xlabel={$t$}, ylabel={$X_t$},
  xmin=1e-6, xmax=1, ymin=-2.1, ymax=2.1,
  ytick={-2,-1,0,1,2},
  tick label style={font=\tiny},
  label style={font=\scriptsize},
  title style={font=\small},
  grid=major, major grid style={gray!20},
]
\addplot[dashed, gray, domain=1e-6:1, samples=2, forget plot] {1};
\addplot[dashed, gray, domain=1e-6:1, samples=2, forget plot] {-1};
\addplot[line width=0.35pt, blue!75!black] table[x=t,y=x1] {Python/data/paths_g075.csv};
\addplot[line width=0.35pt, blue!55!black] table[x=t,y=x2] {Python/data/paths_g075.csv};
\addplot[line width=0.35pt, blue!35!black] table[x=t,y=x3] {Python/data/paths_g075.csv};
\addplot[line width=0.35pt, blue!75!black, densely dotted] table[x=t,y=x4] {Python/data/paths_g075.csv};
\addplot[line width=0.35pt, blue!55!black, densely dotted] table[x=t,y=x5] {Python/data/paths_g075.csv};
\addplot[line width=0.35pt, blue!35!black, densely dotted] table[x=t,y=x6] {Python/data/paths_g075.csv};
\end{semilogxaxis}
\end{tikzpicture}
\end{minipage}\hfill
\begin{minipage}[t]{0.305\textwidth}
\centering
\begin{tikzpicture}[baseline=(current bounding box.north)]
\begin{semilogxaxis}[
  width=\linewidth, height=0.78\linewidth,
  title={(c) switched, $\gamma=1$},
  xlabel={$t$}, ylabel={$X_t$},
  xmin=1e-6, xmax=1, ymin=-2.1, ymax=2.1,
  ytick={-2,-1,0,1,2},
  tick label style={font=\tiny},
  label style={font=\scriptsize},
  title style={font=\small},
  grid=major, major grid style={gray!20},
]
\addplot[dashed, gray, domain=1e-6:1, samples=2, forget plot] {1};
\addplot[dashed, gray, domain=1e-6:1, samples=2, forget plot] {-1};
\addplot[line width=0.35pt, orange!85!black] table[x=t,y=x1] {Python/data/paths_g1.csv};
\addplot[line width=0.35pt, orange!70!black] table[x=t,y=x2] {Python/data/paths_g1.csv};
\addplot[line width=0.35pt, orange!55!black] table[x=t,y=x3] {Python/data/paths_g1.csv};
\addplot[line width=0.35pt, orange!85!black, densely dotted] table[x=t,y=x4] {Python/data/paths_g1.csv};
\addplot[line width=0.35pt, orange!70!black, densely dotted] table[x=t,y=x5] {Python/data/paths_g1.csv};
\addplot[line width=0.35pt, orange!55!black, densely dotted] table[x=t,y=x6] {Python/data/paths_g1.csv};
\end{semilogxaxis}
\end{tikzpicture}
\end{minipage}

\vspace{2pt}

\begin{minipage}[t]{0.328\textwidth}
\centering
\begin{tikzpicture}[baseline=(current bounding box.north)]
\begin{loglogaxis}[
  width=0.90\linewidth, height=0.66\linewidth,
  title={(d) distance to the atom},
  xmin=1e-6, xmax=1, ymin=1e-4, ymax=2,
  xlabel={$t$}, ylabel={$|X_t-y_I|$},
  ytick={1e-3,1e-2,1e-1,1},
  tick label style={font=\tiny},
  label style={font=\scriptsize},
  title style={font=\small},
  grid=major, major grid style={gray!20},
]
\addplot[thick, blue!75!black] table[x=t, y=g025] {Python/data/rates.csv};
\addplot[thick, teal!80!black] table[x=t, y=g075] {Python/data/rates.csv};
\addplot[thick, orange!85!black] table[x=t, y=crit] {Python/data/rates.csv};
\addplot[dashed, black, domain=1e-6:1, samples=2, forget plot] {0.02};
\node[font=\tiny, anchor=south, fill=white, inner sep=1pt]
  at (axis cs:3e-4,0.024) {critical: $2\kappa$};
\node[font=\tiny, anchor=north, fill=white, inner sep=1pt]
  at (axis cs:3e-4,9e-4) {$\gamma=0.25,\,0.75$};
\end{loglogaxis}
\end{tikzpicture}
\end{minipage}\hfill
\begin{minipage}[t]{0.328\textwidth}
\centering
\begin{tikzpicture}[baseline=(current bounding box.north)]
\begin{axis}[
  width=0.86\linewidth, height=0.62\linewidth,
  title={(e) endpoint frequencies},
  xmin=0.4, xmax=3.6, ymin=0.485, ymax=0.545,
  xtick={1,2,3}, xticklabels={$0.25$,$0.5$,$0.75$},
  ytick={0.50,0.51,0.52,0.53},
  xlabel={$\gamma$}, ylabel={$\hat q_+$},
  tick label style={font=\tiny},
  label style={font=\scriptsize},
  title style={font=\small},
  grid=major, major grid style={gray!20},
]
\addplot[only marks, mark=*, mark size=2pt, blue!75!black,
  error bars/.cd, y dir=both, y explicit, error bar style={blue!75!black}]
  table[x=idx, y=qhat, y error plus=err_plus, y error minus=err_minus]
  {Python/data/endpoint_bars.csv};
\addplot[dashed, black, domain=0.4:3.6, samples=2, forget plot] {0.5};
\end{axis}
\end{tikzpicture}
\end{minipage}\hfill
\begin{minipage}[t]{0.328\textwidth}
\centering
\begin{tikzpicture}[baseline=(current bounding box.north)]
\begin{axis}[
  width=0.86\linewidth, height=0.62\linewidth,
  title={(f) terminal offset},
  xmin=-0.008, xmax=0.042, ymin=0, ymax=380,
  xtick={0,0.01,0.02,0.03,0.04},
  ytick={0,100,200,300},
  xticklabel style={/pgf/number format/fixed,
                    /pgf/number format/precision=2,
                    /pgf/number format/zerofill,
                    font=\tiny},
  yticklabel style={/pgf/number format/fixed, font=\tiny},
  xlabel={$X_{t_{\mathrm{stop}}}-y_I$}, ylabel={density},
  tick label style={font=\tiny},
  label style={font=\scriptsize},
  title style={font=\small},
]
\addplot[thick, blue!75!black] table[x=x, y=g075] {Python/data/offset_kde.csv};
\addplot[thick, orange!85!black] table[x=x, y=g1s] {Python/data/offset_kde.csv};
\draw[dotted, thick] (axis cs:0.0017,0) -- (axis cs:0.0017,350);
\draw[dotted, thick] (axis cs:0.02,0) -- (axis cs:0.02,350);
\node[font=\tiny, anchor=south, fill=white, inner sep=1pt]
  at (axis cs:0.0017,290) {$\gamma=0.75$};
\node[font=\tiny, anchor=south, fill=white, inner sep=1pt]
  at (axis cs:0.02,335) {critical};
\end{axis}
\end{tikzpicture}
\end{minipage}
\caption{Simulations of \eqref{eq:main-sde} for $S=\{-1,+1\}$, $d=1$, $T=1$
and initial law $v_T=w(T)$; seed $20261028$, $N=10^4$, $t_{\mathrm{stop}}=10^{-6}$.
(a)--(c) Six sample paths for $\gamma=0.25$, $\gamma=0.75$, and the switched
critical perturbation of Proposition~\ref{prop:critical} ($\kappa=0.01$,
$b=0.05$); dashed lines mark $\pm1$. (d) Median distance $|X_t-y_I|$ to the
selected atom on a logarithmic scale: the two subcritical errors keep decaying,
whereas the critical perturbation levels off at $2\kappa=0.02$ (dashed line).
(e) Endpoint frequencies $\hat q_+$ with Wilson $95\%$ intervals; the dashed line
is $1/2$. (f) Kernel density of the terminal offset $X_{t_{\mathrm{stop}}}-y_I$
for $\gamma=0.75$, centered at $2Mt_{\mathrm{stop}}^{1-\gamma}/\gamma=0.0017$, and
for the switched perturbation, centered at $2\kappa=0.02$ (dotted lines).}
\label{fig:response-energy-illustration}
\end{figure}
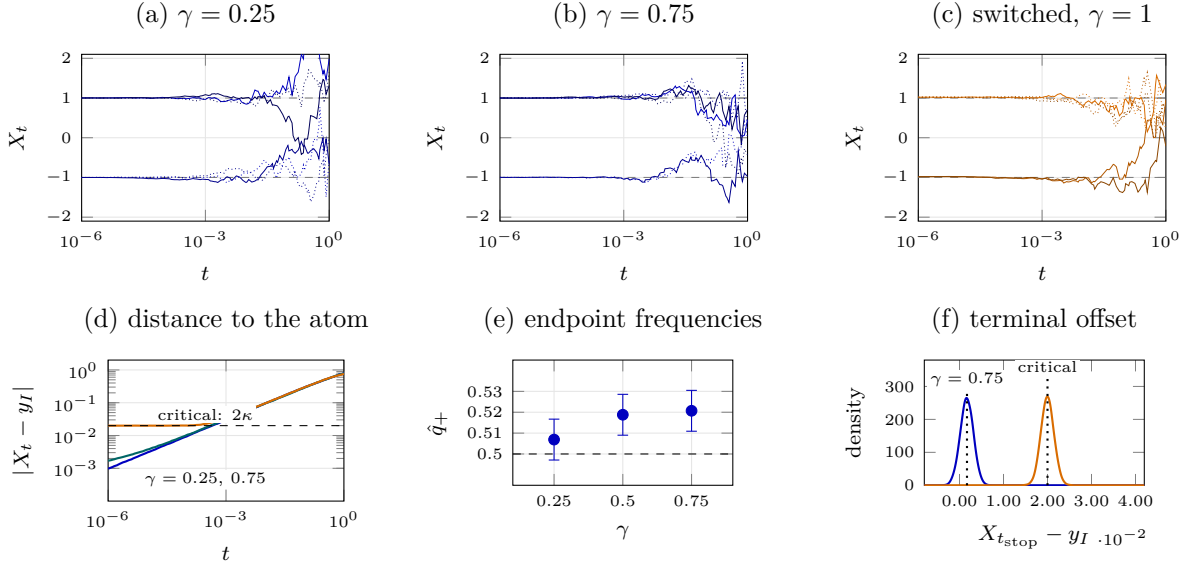

Panels (a)--(c) of Figure~\ref{fig:response-energy-illustration} show six sample
paths in each regime. In all three regimes the trajectory is attracted to a single
atom well before $t=10^{-2}$. The difference between the regimes is not visible on the
scale of the paths but in the terminal offset, which is resolved in panels (d)
and (f). Panel (e) gives the endpoint frequencies: $\hat q_+=0.507$, $0.519$ and $0.521$
for $\gamma=0.25,0.5,0.75$, with Wilson $95\%$ intervals of half-width $0.010$.
The support is $\{-1,+1\}$ in every run, while the weights are shifted modestly
upward by the fixed direction $\xi=+1$. Neither Theorem~\ref{thm:main} nor
Corollary~\ref{cor:positive-terminal-support} predicts the values of the $q_i$,
so this behaviour is compatible with Remark~\ref{rem:fidelity}.

Panel (d) shows the median distance to the selected atom as a function of
$t$ on a logarithmic scale: for the two subcritical errors it keeps decaying,
while for the critical perturbation it levels off at $2\kappa=0.02$, so the
displacement is a limit and not a transient. Panel (f) resolves the terminal
offset. For $\gamma=0.75$ its density is concentrated at $2Mt_{\mathrm{stop}}^{1-\gamma}/\gamma=0.0017$, the leading term
of \eqref{eq:matching}; for the switched perturbation it is concentrated at
$2\kappa=0.020$, in agreement with \eqref{eq:critical-limit}. Both values are
predicted by the asymptotic statements and are not fitted to the simulation.

\begin{table}[t]
\centering
\small
\setlength{\tabcolsep}{4.5pt}
\begin{tabular}{@{}lcccccc@{}}
\hline
$t_{\mathrm{stop}}$ & $\hat q_+(0.25)$ & $\hat q_+(0.5)$ & $\hat q_+(0.75)$
& med$\,|X-y_I|$, $\gamma=0.75$ & med$\,|X-y_I|$, critical & mean $X$ on $+1$ \\
\hline
$10^{-4}$ & 0.507 & 0.519 & 0.521 & $1.03\times10^{-2}$ & $2.02\times10^{-2}$ & 1.0196 \\
$10^{-5}$ & 0.507 & 0.519 & 0.521 & $3.73\times10^{-3}$ & $1.99\times10^{-2}$ & 1.0194 \\
$10^{-6}$ & 0.507 & 0.519 & 0.521 & $1.70\times10^{-3}$ & $2.00\times10^{-2}$ & 1.0197 \\
\hline
\end{tabular}
\caption{Sensitivity to the stopping time, with the parameters of
Figure~\ref{fig:response-energy-illustration}. Columns: endpoint frequencies for
the three persistent errors; median distance to the selected atom for
$\gamma=0.75$ and for the switched perturbation; mean of
$X_{t_{\mathrm{stop}}}$ over the trajectories ending at $+1$ in the switched
case, whose predicted limit is $1+2\kappa=1.02$.}
\label{tab:tstop}
\end{table}

Table~\ref{tab:tstop} reports the dependence on the stopping time. Over three
decades of $t_{\mathrm{stop}}$ the frequencies are unchanged, since the
localization occurs before $t=10^{-4}$. The median distance to the selected atom
for $\gamma=0.75$ decreases from $1.03\times10^{-2}$ to $1.70\times10^{-3}$,
a log--log slope of about $0.39$ over this range. This is intermediate between
the Brownian exponent $1/2$ and the asymptotic exponent $1-\gamma=1/4$ of
Proposition~\ref{prop:matching}: at $t_{\mathrm{stop}}=10^{-6}$ the median
agrees with the leading term $2Mt^{1-\gamma}/\gamma=1.69\times10^{-3}$ to
better than $1\%$, whereas the Brownian scale
$\sqrt{2t_{\mathrm{stop}}}=1.41\times10^{-3}$ differs by about $20\%$. At
$t_{\mathrm{stop}}=10^{-4}$ the two scales themselves differ by a factor
$2.6$, so the experiment illustrates the crossover between them rather than
resolving the asymptotic rate. The critical offset stays at $2\kappa$ within
the Monte Carlo error. Halving
$\varepsilon$ leaves the median critical offset unchanged to $10^{-4}$ and
changes the frequency estimate by less than one Monte Carlo standard error in a
check with $n=2000$ trajectories.

Two caveats. First, all results refer to $t_{\mathrm{stop}}>0$: the offset in
panel (f) is the finite-$t_{\mathrm{stop}}$ value of the asymptotic statements,
and it moves with $t_{\mathrm{stop}}$ as predicted, but no simulation verifies an
almost-sure statement at $t\downarrow0$. Second, the experiment is deliberately
small---one target, one direction $\xi$, and $10^4$ trajectories---and
illustrates the two thresholds and the terminal displacement rather than
benchmarking a sampler. The script \texttt{simulate\_reverse\_sde.py} and the data
files used for Figure~\ref{fig:response-energy-illustration} and
Table~\ref{tab:tstop} are included in the source folder \texttt{Python/}.

\section{Conclusion and discussion}\label{sec:conclusion}

For finite empirical targets, the contraction-weighted condition
$t\int_t^T\frac{E(\sigma)}{\sigma}\dd \sigma \to0$ yields almost-sure convergence to a target atom and an explicit
trajectory bound throughout a deterministic terminal window $0<t\le b$, near the same training sample. An exact space
and time change and a global Lyapunov identity establish convergence; local
propagation then yields the uniform terminal estimate. The weighted condition
is strictly weaker than $t\int_t^T E(\sigma)^2\dd \sigma \to0$ and includes sparse high-amplitude errors
excluded by that squared-error condition. Both localization and full
terminal support hold for arbitrary initial probability laws, without entropy
or moment bounds. Every atom has positive terminal probability,
and $\operatorname{supp}(v_0^*)=S$. Localization can persist
with infinite error energy and mutually singular perturbed and exact path laws.

For persistent power-law errors, the critical scale is $t^{-1}$.
For any initialization probability law, Proposition~\ref{prop:critical} constructs
perturbations with arbitrarily small positive amplitude whose terminal limit
lies at distance $2\kappa$ from the training set with probability at least
$1-\alpha$. The switching window depends on $v_T$, $\kappa$ and $\alpha$, so these
conclusions concern the corresponding chosen process. Together, the results
show how terminal contraction can suppress divergent score errors, while
critical perturbations can alter the generated endpoint.

For diffusion modeling in this finite empirical setting, the infinite-energy
examples show that a divergent terminal score error does not by itself prevent
convergence to the training set under the stated error conditions. Its growth
alone therefore does not establish that a continuous-time sampler has escaped
terminal attraction to training points. The denoiser correction
$\delta_r^\mu=h_{T-r}e_r$, with $h_t\sim2t$, distinguishes score divergence
from the effective correction: a divergent score error can correspond to a
vanishing correction, and Remark~\ref{rem:denoiser} gives a sufficient
criterion through a deterministic correction envelope.

Terminal localization and distributional fidelity are distinct:
the atom probabilities need not remain $1/N$, and their errors are not
estimated by the localization theorem; see~\eqref{eq:terminal-tv}.
Example~\ref{ex:changed-weights} gives an explicit bounded perturbation
that changes those probabilities while preserving the two-point support
and the exact sampler's initialization. The conclusions concern terminal limits of
the continuous-time model and do not directly extend to finite-step discrete
samplers. For qualitative almost-sure localization,
Proposition~\ref{prop:global-localization} and
Corollary~\ref{cor:mean-integrated} now permit realized weighted integrals
and mean-integrated control under the generated law. Transferring
reference-law learning bounds to that law, obtaining quantitative bounds
under these weaker hypotheses, and extending the mechanism to non-atomic
targets remain directions for future work.

\appendix
\renewcommand{\thesection}{\Alph{section}}
\section{Well-posedness and terminal path-law comparison}\label{app:prelim}

We record the well-posedness used on truncated intervals and for prescribed
perturbations, then compare the full terminal path laws for an explicit
infinite-energy perturbation.

\subsection{Well-posedness before the terminal time}\label{app:wellposed}

Fix $0<\delta<T$ and put $H=T-\delta$. The finite Gaussian mixture $w(\cdot,t)$
has positive component variances for $t>0$, so its score is smooth in the
spatial variable. Recall from \eqref{eq:score} that
$\mu(x,t)=\sum_{i=1}^N \omega_i(x,t)a_ty_i$ is the conditional mean of the moving
centers, and $s(x,t)=\nabla_x\log w(x,t)=(\mu(x,t)-x)/h_t$.
Here $\|A\|_{\Op}=\sup_{|z|=1}|Az|$ is the Euclidean-induced operator norm.
Direct differentiation gives
\[
\nabla_x\mu(x,t)=h_t^{-1}\Sigma(x,t),\qquad
\Sigma(x,t)=a_t^2\sum_{i=1}^N \omega_i(x,t)y_iy_i^\top
-\mu(x,t)\mu(x,t)^\top.
\]
The matrix $\Sigma(x,t)$ is the conditional covariance of the moving centers
given the noisy state $x$; since $|a_ty_i|\le R$, it satisfies
$0\preceq\Sigma(x,t)\preceq R^2I_d$ and $|\mu(x,t)|\le R$. Consequently,
for $0\le r\le H$ and $t=T-r$,
\begin{equation}\label{eq:lipschitz}
\begin{aligned}
\|\nabla_x[x+2s(x,t)]\|_{\Op}
&\le L_\delta\triangleq1+2/h_\delta+2R^2/h_\delta^2,\\
|x+2s(x,t)|&\le(1+2/h_\delta)|x|+2R/h_\delta.
\end{aligned}
\end{equation}
Thus the exact drift is uniformly globally Lipschitz with linear growth on
this truncated interval.

For a prescribed progressively measurable perturbation satisfying
\eqref{eq:envelope}, local square integrability of $E$ gives
$\int_0^H|e_\ell|\dd\ell\le\sqrt{H\int_\delta^T E(\sigma)^2\dd\sigma}<\infty$ almost surely.
For every $\F_0$-measurable $\R^d$-valued $\widehat X_0$, the process
$A_r=\widehat X_0+2\int_0^r e_\ell\dd\ell+\sqrt2B_r$ is therefore continuous and adapted.
Pathwise Picard iteration for
$\widehat X_r=A_r+\int_0^r[\widehat X_\ell+2s(\widehat X_\ell,T-\ell)]\dd\ell$, together with
\eqref{eq:lipschitz} and Gronwall's inequality, gives a pathwise unique strong
solution on $[0,H]$ with continuous adapted paths, without a moment assumption
on $\widehat X_0$. Pathwise uniqueness makes these solutions agree on overlapping
intervals for the same perturbation, initial value and Brownian motion;
a countable exhaustion of $[0,T)$ gives a unique continuous adapted solution
there. For deterministic perturbations and each fixed initialization law, uniqueness
in law holds on every truncated interval, as used below. For general feedback errors, existence
remains part of Assumption~\ref{ass:localization}\textup{(\ref{ass:local})}.

\subsection{Local equivalence and terminal path-law singularity}

Equivalence of path laws on every interval strictly before the terminal
generation time does not imply equivalence on the sigma-field generated by
the entire trajectory. The following example has mutually singular full path
laws and still satisfies the localization criterion. We define the density
process on a common uncompleted canonical path space.

Work on $\Omega_0\triangleq C([0,T);\R^d)$ with the increasing-time coordinate process
$\widehat X_r(\omega)\triangleq\omega(r)$ and
canonical filtration without completion,
$\F_r^0\triangleq \sigma(\widehat X_\ell:0\le\ell\le r)$.
On this space set $X_t:=\widehat X_{T-t}$ for $0<t\le T$; the canonical filtration
continues to be indexed by $r$.
The sigma-field generated by the whole trajectory before $r=T$ is
$\F_{T-}^0\triangleq \sigma(\bigcup_{r<T}\F_r^0)$.
The countable-coordinate construction in the proof of
Proposition~\ref{prop:global-localization} defines an $\F_{T-}^0$-measurable
label here, agreeing with the terminal atom index whenever the coordinate
path converges to a point of $S$. Fix the deterministic perturbation
\begin{equation}\label{eq:nonui-error}
    e_r\triangleq M(T-r)^{-\gamma}\xi,\qquad M>0,\quad|\xi|=1,\quad
    \frac12\le\gamma<1.
\end{equation}
Let $P$ and $Q$ be the exact and perturbed laws on this common canonical
space, both initialized from $w(T)$, with perturbation \eqref{eq:nonui-error}
under $Q$.
These laws exist by the well-posedness argument in Appendix~\ref{app:wellposed}. The reconstructed noise
\[
B_r^P\triangleq \frac1{\sqrt2}
\left(\widehat X_r-\widehat X_0-\int_0^r[\widehat X_\ell+2s(\widehat X_\ell,T-\ell)]\dd\ell\right)
\]
is Brownian under $P$. Since $e$ is $C^1$ on each $[0,r]$, $r<T$, we fix
a version of $\int_0^r e_\ell\cdot\dd B_\ell^P$ by integration by parts as
$e_r\cdot B_r^P-\int_0^r\dot e_\ell\cdot B_\ell^P\dd\ell$, using $B_0^P=0$.
This expression is measurable with respect to the uncompleted canonical
filtration. It defines the same continuous exponential process below as a
function of the coordinate path under both laws.

\begin{lemma}[Mutual singularity despite terminal localization]\label{lem:nonui}
    For every $r<T$, the restrictions of $P$ and $Q$ to $\F_r^0$ are equivalent.
    Their Radon--Nikodym density is given by
    \begin{equation}\label{eq:Ldef}
        \frac{\dd Q|_{\F_r^0}}{\dd P|_{\F_r^0}}=L_r,\qquad
        L_r\triangleq\exp\left\{\sqrt2\int_0^r e_\ell\cdot\dd B_\ell^P
        -\int_0^r|e_\ell|^2\dd\ell\right\}.
    \end{equation}
    The process $(L_r)_{r<T}$ is a positive $P$-martingale with $\E_PL_r=1$.
    With
    \begin{equation}\label{eq:Hdef}
        H_r\triangleq 2\int_0^r|e_\ell|^2\dd\ell
        =\begin{cases}
            2M^2\log\dfrac{T}{T-r},&\gamma=\frac12,\\[4pt]
            \dfrac{2M^2}{2\gamma-1}
            \bigl[(T-r)^{1-2\gamma}-T^{1-2\gamma}\bigr],
            &\frac12<\gamma<1,
        \end{cases}
    \end{equation}
    one has $H_r\uparrow\infty$ as $r\uparrow T$ and
    \[
    L_r\longrightarrow0\quad P\text{-a.s.},
    \qquad
    L_r\longrightarrow\infty\quad Q\text{-a.s.}
    \]
    Moreover, $P\perp Q$ on $\F_{T-}^0$. The process $(L_r)_{r<T}$ is not
    uniformly integrable under $P$, and no terminal Radon--Nikodym density
    $\dd Q/\dd P$ exists on $\F_{T-}^0$.
\end{lemma}

\begin{proof}[Proof of Lemma~\ref{lem:nonui}]
\stepheading{Equivalence before the terminal time.}
For each $H<T$, the deterministic energy $\int_0^H|e_\ell|^2\dd\ell$ is finite.
Novikov's condition makes $L$ a positive $P$-martingale on $[0,H]$.
Under $L_H\cdot P|_{\F_H^0}$,
\[
 B_r^Q=B_r^P-\sqrt2\int_0^r e_\ell\dd\ell
\]
is Brownian and $\widehat X$ solves the perturbed equation. The conditional identity
$\E_P[L_H\mid\F_0^0]=1$ preserves the initialization law. Uniqueness in law
therefore identifies this measure with $Q|_{\F_H^0}$. Hence the restrictions
are equivalent and have density $L_r$ for every $r<T$. Compatible continuous
versions along a countable exhaustion $H_m\uparrow T$ make these identities
simultaneous outside a set null under both laws.

\proofstep{Limits of the density process.}
The stochastic integral in $\log L_r$ has quadratic variation $H_r$.
Dambis--Dubins--Schwarz and the Brownian strong law give
\[
 \frac{\log L_r}{H_r}
 =\frac{W^P_{H_r}}{H_r}-\frac12\longrightarrow-\frac12
 \qquad P\text{-a.s.}
\]
Under $Q$, the common pathwise definition of the integral gives
\begin{equation}\label{eq:Ltwo}
 \log L_r=\sqrt2\int_0^r e_\ell\cdot\dd B_\ell^Q+\frac12H_r,
\end{equation}
so the same argument yields
\[
 \frac{\log L_r}{H_r}
 =\frac{W^Q_{H_r}}{H_r}+\frac12\longrightarrow\frac12
 \qquad Q\text{-a.s.}
\]
Since $H_r\uparrow\infty$, these imply the asserted limits of $L_r$.

\proofstep{A separating event.}
For $r_m=T(1-2^{-m})$, define
\[
 A_0=\{\lim_{m\to\infty}L_{r_m}=0\}
 =\bigcap_{k\ge1}\bigcup_{j\ge1}\bigcap_{m\ge j}
 \{L_{r_m}<1/k\}\in\F_{T-}^0.
\]
The preceding limits give $P(A_0)=1$ and $Q(A_0)=0$.
Thus $P\perp Q$ on $\F_{T-}^0$, and no terminal density $\dd Q/\dd P$
exists there. Finally, $L_r\to0$ $P$-almost surely while $\E_PL_r=1$,
so the density martingale is not uniformly integrable.
\end{proof}


\end{document}